[1994/06/01]% LaTeX date must June 1994 or later

\documentclass[12pt]{amsart}

\usepackage{amsmath,amssymb,enumerate}
\usepackage{epsfig,fancyhdr,color}
\usepackage{amssymb}
\usepackage{amsmath,amsthm}
\usepackage{latexsym}
\usepackage{amscd}
\usepackage{psfrag}
\usepackage{graphicx}
\usepackage[latin1]{inputenc}
\usepackage[all]{xy}
\usepackage[mathcal]{eucal}

        [1994/05/12 v2.2b Standard AMS font definitions]
\DeclareFontFamily{U}{msa}{}
\DeclareFontShape{U}{msa}{m}{n}
  { <5> <6> <7> <8> <9> gen * msam
    <10> <10.95> <12> <14.4> <17.28> <20.74> <24.88> msam10}{}
        [1994/05/12 v2.2b Standard AMS font definitions]
\DeclareFontFamily{U}{msb}{}
\DeclareFontShape{U}{msb}{m}{n}
  { <5> <6> <7> <8> <9> gen * msbm
    <10> <10.95> <12> <14.4> <17.28> <20.74> <24.88> msbm10}{}

\newtheorem{cor}[subsection]{Corollary}
\newtheorem{lem}[subsection]{Lemma}
\newtheorem{prop}[subsection]{Proposition}
\newtheorem{rem}[subsection]{Remark}
\newtheorem{thm}[subsection]{Theorem}
\newtheorem{defn}[subsection]{Definition}
\newtheorem{exam}[subsection]{Example}

\theoremstyle{definition}

\numberwithin{equation}{section}

\title[The Complex of Domains]{Automorphisms of the Complex of Domains}

\author[J. D. McCarthy]{John D. McCarthy}
\address{Department of Mathematics, Michigan State University, East Lansing, MI 48824-1027, USA and Max-Plank-Institut f\"ur Mathematik, Vivatsgasse 7, 53111 Bonn, Germany}
\email{mccarthy@math.msu.edu}
\thanks{The authors thank the Max Planck Institute for Mathematics, Bonn for the excellent conditions provided for their stay at this institution, during which this paper was written.}

\author[A. Papadopoulos]{Athanase Papadopoulos}
\address{Institut de Recherche Math\'{e}matique Avanc\'{e}e, Universit\'{e} Louis Pasteur et CNRS, 7 rue Ren\'e-Descartes,  67084 Strasbourg Cedex, France and Max-Plank-Institut f\"ur Mathematik, Vivatsgasse 7, 53111 Bonn, Germany}
\email{papadopoulos@math.u-strasbg.fr}

\date{\today}


\newcommand{\abstracttext}{In this paper, we study a flag complex which is naturally associated to the Thurston theory of surface diffeomorphisms for compact connected orientable surfaces with boundary \cite{mccpap}. The various pieces of the Thurston decomposition of a surface diffeomorphism, thick domains and annular or thin domains, fit into this flag complex, which we call the complex of domains. The main result of this paper is a computation of the group of automorphisms of this complex. Unlike the complex of curves, introduced by Harvey \cite{harvey}, for which, for all but a finite number of exceptional surfaces, by the works of Ivanov \cite{ivanov1}, Korkmaz \cite{korkmaz}, and Luo \cite{luo}, all automorphisms are geometric (i.e. induced by homeomorphisms), the complex of domains has nongeometric automorphisms, provided the surface in question has at least two boundary components. These nongeometric automorphisms of the complex of domains are associated to certain edges of the complex which are naturally associated to biperipheral pairs of pants on the surface in question. We project the complex of domains to a natural subcomplex of the complex of domains by collapsing each biperipheral edge onto the unique vertex of that edge which is represented by a regular neighborhood of a biperipheral curve and, thereby, reduce the computation in question to computing the group of automorphisms of this subcomplex, which we call the truncated complex of domains. Finally, we prove that the group of automorphisms of the truncated complex of domains is the extended mapping class group of the surface in question and, obtain, thereby, a complete description of the group of automorphisms of the complex of domains.}

\begin{document}

\begin{abstract}  \abstracttext  \end{abstract}

\maketitle

\section{Introduction}

Let $S = S_{g,b}$ be a connected compact orientable surface of genus $g$ with $b$ boundary components. Let $\partial S$ denote the boundary of $S$. The {\it mapping class group of $S$},\index{mapping class group} $\Gamma ( = \Gamma_{g,b} = \Gamma(S))$, is the group of isotopy classes of orientation-preserving self-homeomorphisms of $S$. The {\it extended mapping class group of $S$},\index{extended mapping class group} \index{mapping class group!extended} $\Gamma^*$, is the group of isotopy classes of self-homeomorphisms of $S$. Note that $\Gamma$ is a normal subgroup of index $2$ in $\Gamma^*$. 

The study of these groups, $\Gamma$ and $\Gamma^*$. has used their action on various abstract simplicial complexes, each of which encodes combinatorial information about the relationship which certain subspaces of $S$ bear to one another.  For instance, the curve complex, $C(S)$, which was introduced by W. Harvey \cite{harvey}, captures the combinatorial complexity of the set of isotopy classes of essential unoriented simple closed curves on $S$. 

In recent joint work, we have begun the study of a new complex on which $\Gamma^*$ acts \cite{mccpap}. This complex is naturally associated to the Thurston theory of surface diffeomorphisms for compact connected orientable surfaces with boundary. The various pieces of the Thurston decomposition of a surface diffeomorphism, thick domains and annular or thin domains, fit into this flag complex, which we call the complex of domains. 

More precisely, a {\it domain on $S$} is a nonempty connected compact embedded surface in $S$ which is not equal to $S$ and each of whose boundary components is either contained in $\partial S$ or is essential on $S$.  The vertex set $D_0(S)$ of $D(S)$ is the set of isotopy classes of  domains on $S$. An $n$-simplex of $D(S)$ is a set of $n + 1$ distinct vertices of $D(S)$ which can be represented by disjoint domains of $S$. 

The main result of this paper is our computation of the group of automorphisms of $D(S)$. Unlike the celebrated complex of curves introduced by Harvey \cite{harvey}, for which, for all but a finite number of exceptional surfaces, by the works of Ivanov \cite{ivanov1}, Korkmaz \cite{korkmaz}, and Luo \cite{luo}, all automorphisms are geometric (i.e. induced by homeomorphisms), the complex of domains has nongeometric automorphisms, provided $S$ has at least two boundary components. These nongeometric automorphisms of $D(S)$ are associated to natural biperipheral edges of $D(S)$. 

More precisely, a {\it biperipheral edge of $D(S)$} is an edge of $D(S)$ whose vertices are represented by a {\it biperipheral pair of pants $X$ on $S$}, (i.e. a domain $X$ on $S$ which is homeomorphic to a sphere with three holes having exactly two of its boundary components in $\partial S$) and a regular neighborhood $Y$ of the remaining boundary component of $X$ on $S$, the unique essential boundary component of $X$ on $S$. The corresponding nongeometric automorphism of $D(S)$, which we call a {\it simple exchange of $D(S)$} exchanges the two vertices of $D(S)$ corresponding to $X$ and $Y$ and fixes every other vertex of $D(S)$.  

We prove for most surfaces that every automorphism of $D(S)$ preserves the set $\mathcal{E}$ of all biperipheral edges of $D(S)$ and, hence, induces an automorphism of the subcomplex of $D(S)$ which is obtained from $D(S)$ by removing each vertex of $D(S)$ corresponding to a biperipheral pair of pants on $S$. We call this subcomplex the {\it truncated complex of domains on $S$} and we denote it by $D^2(S)$. In this way, we obtain a natural homomorphism $\rho : Aut(D(S)) \rightarrow Aut(D^2(S))$. 

Studying this homomorphism, $\rho : Aut(D(S)) \rightarrow Aut(D^2(S))$, we prove that it is surjective and that its kernel, which we call {\it the Boolean subgroup $B_\mathcal{E}$ of $D(S)$}, consists of involutions $\varphi_{\mathcal{F}} : D(S) \rightarrow D(S)$, defined for each subset $\mathcal{F}$ of $\mathcal{E}$, which interchange the two vertices of each edge of $D(S)$ in $\mathcal{F}$, and fix every vertex of $D(S)$ which is not a vertex of an edge of $D(S)$ in $\mathcal{F}$. In this way, we see that $B_\mathcal{E}$ is naturally isomorphic to the Boolean algebra $\mathcal{B}(\mathcal{E})$ of all subsets $\mathcal{F}$ of $\mathcal{E}$ and, thereby, exhibit $Aut(D(S))$ as an extension of $Aut(D^2(S))$ by the Boolean algebra $\mathcal{B}(\mathcal{E})$. 

Studying $Aut(D^2(S))$, we prove that the natural representation $\rho : \Gamma^*(S) \rightarrow Aut(D^2(S))$, arising by induction from the natural action of $\Gamma^*(S)$ on $D(S)$, is an isomorphism, completing our computation of $Aut(D(S))$, expressed as follows in our main result. 

\begin{thm}Suppose that $S$ is not a sphere with at most four holes, a torus with at most two holes, or a closed surface of genus two. Then we have a natural commutative diagram of exact sequences: 
 \[ \begin{array}{cccccccccc}
    1 &\longrightarrow & \mathcal{B}(\mathcal{E})   & \longrightarrow &
    \mathcal{B}(\mathcal{E}) \rtimes \Gamma^*(S) & \longrightarrow  &  
    \Gamma^*(S) & 
    \longrightarrow & 1\\
   \   & \  & \simeq\big\downarrow & \  & \simeq\big\downarrow & \  & 
   \simeq\big\downarrow
      & \ & \  \\
     1 & \longrightarrow & B_\mathcal{E}    & \longrightarrow &
    Aut(D(S)) & \longrightarrow  & Aut(D^2(S)) & \longrightarrow & 1 
     \end{array} \]
    \end{thm}
    
Note that, in the case where the surface $S$ has at most one boundary component, we have $D^2(S) = D(S)$ and, therefore, we have a natural isomorphism $\Gamma^*(S) \stackrel{\simeq}{\rightarrow} Aut(D(S))$.

The plan of the rest of this paper is as follows.

  Section 2 contains some basic pronciples on surfaces, subsurfaces, curves,  and geometric intersection numbers that will be used in the later sections of this paper. We have included some of the proofs because they seem to be nonexistent in the literature. Of course, a reader who is an expert in the theory of surfaces will skip these proofs.

  Section 3 contains a short introduction to abstract simplicial complexes, in which we define some basic  simplicial notions that we use in our paper.

In Section 4, we introduce the notion of an exchangeable pair of vertices and of an exchange automorphism of a simplicial complex $K$. We give necessary and sufficient conditionf for a pair of vertices to be exchangeable. We define certain special subgroups of the automorphism group of $K$, which we call Boolean subgroups. Such a group is isomorphic to the Boolean algebra of a collection of  subsets of $K$ consisting of exchangeable pairs of vertices. The exchange automorphisms and the Boolean subgroups will be used in an essential way in the results on the automorphism group of the complex of domains, that we obtain in the sequel.

In Section 5, we introduce the main three simplicial complexes associated to a surface $S$ that we study in this paper: the complex of curves, $C(S)$, the complex of domains, $D(S)$, and the truncated complex of domains, $D^2(S)$. There are natural inclusion maps $C(S)\hookrightarrow D(S)$ and $C(S)\hookrightarrow D^2(S)$ and a natural projection map $D(S)\to D^2(S)$.

In Section 6, we give a simplicial characterization of annular vertices in $D^2(S)$ and we obtain as a result that in the case where the surface $S$ is not a torus with one hole, every simplicial automorphism of $D^2(S)$ restricts to a simplicial automorphism of the subcomplex of $D^2(S)$ that is naturally identified with $C(S)$.

In Section 7, we intruduce the annular link  $Ann(x)$ of a vertex $x$ of $D(S)$. This is the subcomplex  consisting of the simplices in the link of $x$ all of whose vertices are annular. The main result of this section is, provided the surface $S$ is not a sphere with four holes nor a torus with one hole, a characterization of the equality $Ann(x)=Ann(y)$, for vertices $x$ and $y$ of $D(S)$, in terms of the existence of a pair $(X,Y)$ of domains satisfying a certain topological property, and  representing the pair of vertices $(x,y)$.   Thus, while the result of the previous section can be considered as translating  topological information on the surface $S$ into simplicial data, the result of the present section can be considered as translating simplicial information into topological data. This back and forth translation between topological and simplicial information is at the heart of the work done in this paper.

In Section 8, we use the preceding results to prove that in the case where the surface $S$ is not a sphere with four holes or a torus with at most two holes or a closed surface of genus two, all automorphisms of $D^2(S)$ are geometric.

In Section 9, we obtain a simplicial characterization of biperipheral edges in $D(S)$. One of the results that we obtain says that if we exclude a certain finite number of cases that we analyse completely in the paper \cite{mccpap}, the two vertices of a biperipheral edge in $D(S)$ are characterized by the fact that they have the same stars.

In Section 10, we apply the results of the previous section, together with those of Section 4, to study exchange automorphisms of $D(S)$.

Section 11 contains the main results of this paper, that include a complete description of the automorphism group of $D(S)$, for all surfaces $S$ except a finite number of special surfaces which, as mentioned above, we study separately in the paper \cite{mccpap}.

\section{Preliminaries on surfaces} 

For $g \geq 0$ and $b \geq 0$, let $S = S_{g,b}$ be a connected compact orientable surface of genus $g$ with $b$ boundary components. We say that {\it $S$ is a surface of genus $g$ with $b$ holes}. Note that $S_{g,b}$ is a closed surface of genus $g$ if and only if $b = 0$; $S_{0,1}$ is a disc; $S_{0,2}$ is an annulus; $S_{0,3}$ is a pair of pants; $S_{0,b}$ is a sphere with $b$ holes; and $S_{1,b}$ is a torus with $b$ holes.  

Let $\partial S$ denote the boundary of $S$. For convenience, we shall index the $b$ boundary components of $S$, $\partial_i$, $1 \leq i \leq b$.

For each collection $C$ of subsets of $S$, the {\it support of $C$ on $S$} is the union in $S$ of the subsets of $S$ in the collection $C$. We denote the support of $C$ on $S$ by $|C|$.

An isotopy from a homeomorphism $H : S \rightarrow S$ of $S$ to a homeomorphism $H' : S \rightarrow S$ of $S$ is a map $\varphi : S \times [O,1] \rightarrow S$ such that the maps $\varphi_t : S \rightarrow S$, $0 \leq t \leq 1$, are homeomorphisms of $S$, $\varphi_0 = H : S \rightarrow S$, and $\varphi_1 = H' : S \rightarrow S$.

Throughout this paper, all isotopies between subspaces of $S$ will be ambient isotopies. More precisely, if $X$ and $Y$ are subspaces of $S$, an {\it isotopy from $X$ to $Y$} is an isotopy $\varphi : S \times [0,1] \rightarrow S$ from the identity map $\varphi_0 = id_S : S \rightarrow S$ of $S$ to a homeomorphism $\varphi_1 : S \rightarrow S$ of $S$ such that $\varphi_1(X) = Y$.

We denote the isotopy class of a homeomorphism $H : S \rightarrow S$ of $S$ by $[H]$ and the isotopy class of a subspace $X$ of $S$ by $[X]$. 

A {\it curve} on $S$ is an embedded connected closed one-dimensional submanifold of the interior of $S$. 

Let $\alpha$ be a curve on $S$. We say that {\it $\alpha$ is $k$-peripheral on $S$} if there exists a sphere with $k + 1$ holes, $X$, on $S$ such that $\alpha$ is a boundary component of $X$ and every other boundary component of $X$ is a boundary component of $S$.  We say that {\it $\alpha$ is essential on $S$} if it is neither $0$-peripheral nor $1$-peripheral on $S$. In other words, $\alpha$ is essential on $S$ if and only if there does not exist a disc on $S$ whose boundary is equal to $\alpha$ or an annulus $A$ on $S$ whose boundary is equal to the union of $\alpha$ with a boundary component of $S$. 

If $S$ is a sphere with at most three holes, then there are no essential curves on $S$. Otherwise, there are infinitely many pairwise nonisotopic essential curves on $S$. 

Suppose that $\alpha$ and $\beta$ are disjoint essential curves on $S$. Then $\alpha$ is isotopic to $\beta$ on $S$ if and only if there exists an annulus $A$ on $S$ such that the boundary of $A$ is equal to $\alpha \cup \beta$. 

\begin{defn} A collection of pairwise disjoint essential curves on $S$ is a {\it system of curves on $S$} if the curves in the collection are pairwise nonisotopic. Note that every subcollection of a system of curves on $S$ is a system of curves on $S$. \label{defn:systemcurves} \end{defn}

Note that every system of curves on a closed torus has exactly one curve. 

Let $C$ be a finite collection of curves on $S$ and $S_C$ be the surface obtained from $S$ by cutting $S$ along $|C|$ \cite{ivanov2}. 

\begin{defn} A {\it pants decomposition of $S$} is a collection of curves $C$ on $S$ such that each component of $S_C$ is a pair of pants (i.e. a sphere with three holes).\label{defn:pantsdecomposition} \end{defn} 

Note that every pants decomposition of $S$ is a system of curves on $S$. 

Note that $S$ has a pants decomposition if and only if $S$ is not a sphere with at most two holes nor a closed torus. Moreover, on such a surface $S$, if $C$ is a system of curves on $S$, then there exists a pants decomposition on $S$ containing $C$.    

A nonempty system of curves $C$ on $S$ is a maximal system of curves on $S$ if and only if $S$ is a closed torus and $C$ consists of exactly one nonseparating curve on $S$ or $S$ is not a closed torus and $C$ is a pants decomposition of $S$, in which case the number of curves in $C$ is equal to $3g - 3 + b$.

Suppose that $C$ is a pants decomposition of $S$. Let $R$ be a regular neighborhood of $|C|$. Note that the closure of the complement of $R$ on $S$ is a disjoint union of pairs of pants on $S$. These pairs of pants are called {\it pairs of pants of $C$}. Note that, with this definition, the pairs of pants of a pants decomposition are defined only up to isotopy relative to $C$.

Suppose that $P$ is a pair of pants of a pants decomposition $C$ of $S$. Note that $P$ is contained in a unique component $Q$ of $S_C$. We say that {\it $P$ is an embedded pair of pants of $C$} if the natural quotient map $\pi : S_C \rightarrow S$ embeds the pair of pants $Q$ in $S$.  

Let $C$ be a pants decomposition of $S$. We say that {\it $C$ is an embedded pants decomposition of $S$} if every pair of pants of $C$ is embedded. For example, if $S$ is a closed surface of genus two and $C$ is a disjoint union of three nonisotopic nonseparating curves on $S$, then $C$ is an embedded pants decomposition of $S$.

In the rest of this paper, unless otherwise indicated, all curves will be assumed to be essential. 

\begin{defn} Let $\alpha$ and $\beta$ be curves on $S$. The {\it geometric intersection number of $\alpha$ and $\beta$ on $S$} is the minimum number $i_S(\alpha,\beta)$ of points in $\alpha' \cap \beta'$ where $\alpha'$ and $\beta'$ are curves on $S$ which are isotopic to $\alpha$ and $\beta$ on $S$. \label{defn:geointnum} \end{defn}

\begin{defn} Let $\{\alpha_1,...,\alpha_n\}$ be a collection of curves on $S$. We say that {\it $\{\alpha_1,...,\alpha_n\}$ is in minimal position on $S$} if the geometric intersection number of $\alpha_i$ and $\alpha_j$ on $S$ is equal to the number of points in $\alpha_i \cap \alpha_j$, $1 \leq i < j \leq n$.\label{defn:minpos} \end{defn}

We recall that if $\{\alpha_1,...,\alpha_n\}$ is a finite collection of curves on $S$, then there exists a collection $\{\beta_1,...,\beta_n\}$ of curves on $S$ such that $\alpha_i$ is isotopic to $\beta_i$ on $S$, $1 \leq i \leq n$ and $\{\beta_1,...,\beta_n\}$ is in minimal position on $S$. A proof of this fact can be obtained by using an innermost disk argument. Alternatively, one can prove this fact by using hyperbolic geometry \cite{flp}.

\begin{prop} Let $C$ be a collection of disjoint essential curves on $S$. Let $\alpha \in C$. Then there exists a curve $\gamma$ on $S$ such that for each $\beta \in C$, $i(\beta,\gamma) \neq 0$ if and only if $\alpha$ is isotopic to $\beta$ on $S$. 
\label{prop:detectisotopy} \end{prop}

\begin{proof} Since $\alpha$ is an essential curve on $S$, $S$ is not a sphere with at most three holes. 

Suppose, on the one hand, that $S$ is a closed torus. Since any two disjoint essential curves on $S$ are isotopic, $\alpha$ is isotopic to $\beta$ on $S$ for each $\beta$ in $C$. Since any essential curve on $S$ is nonseparating on $S$, $\alpha$ is a nonseparating curve on $S$. Hence, there exists a curve $\gamma$ on $S$ that intersects $\alpha$ transversely and has exactly one point of intersection with $\alpha$. It follows that $i(\beta,\gamma) = i(\alpha,\gamma) = 1$ for each $\beta$ in $C$.

Suppose, on the other hand, that $S$ is not a closed torus. Then we may choose a pants decomposition $D$ on $S$ such that $\alpha$ is an element of $D$ and each element $\beta$ of $C$ is isotopic to a unique element $\beta'$ of $D$. Since $\alpha$ is a curve of the pants decomposition $D$ of $S$, there exists a curve $\gamma$ on $S$ such that $i(\alpha,\gamma) \neq 0$ and $i(\delta,\gamma) = 0$ for each curve $\delta$ in $D \setminus \{\alpha\}$. Indeed, we can choose $\gamma$ so that $i(\alpha,\gamma) = 1$, if there is a nonembedded pair of pants $P$ of $D$ having two of its boundary components isotopic to $\alpha$, and, otherwise, $i(\alpha,\gamma) = 2$ \cite{flp}. 

It follows that for each $\beta$ in $C$, $i(\beta,\gamma) = i(\beta',\gamma) \neq 0$ if and only if $\alpha$ is isotopic to $\beta$ on $S$ (i.e. if and only if $\alpha = \beta'$).\end{proof}

\begin{prop} Suppose that $S$ is not a sphere with at most three holes. Let $H : S \rightarrow S$ be a homeomorphism of $S$. If $H$ preserves the isotopy class of every essential curve on $S$, then $H$ is 
orientation-preserving. 
\label{prop:orpres} \end{prop}

\begin{proof} First, consider the case where the genus of $S$ is zero. Since $S$ is not a sphere with at most three holes, there exists a sphere with four holes $X$ embedded in $S$ such that three of the four boundary components of $X$ are boundary components $C_1$, $C_2$, and $C_3$ of $S$ and the remaining boundary component $C_0$ of $X$ is either a boundary component of $S$ or an essential curve on $S$. 

As in Section 4.2 of \cite{ivmcc}, we recall the lantern relation discovered by M. Dehn \cite{dehn} and rediscovered and popularized by D. Johnson \cite{johnson}. To do this, we choose an orientation on $S$. Let $\alpha_{ij} = \alpha_{ji} $, $1 \leq i < j \leq 3$ be an arc on $S$ joining $C_i$ to $C_j$. Suppose that $\alpha_{12}$, $\alpha_{23}$, and $\alpha_{31}$ are disjoint. Note that the surface obtained from $X$ by cutting along $\alpha_{12} \cup \alpha_{23} \cup \alpha_{31}$ contains a unique component $D$ which is a disc. Suppose that $D$ is on the left of $\alpha_{12}$ as we travel along $\alpha_{12}$ from $C_1$ to $C_2$, as in Figure \ref{fig:lantern}. Let $C_{ij}=C_{ji}$ be the unique essential boundary component of a regular neighborhood $P_{ij}=P_{ji}$ in $X$ of $C_i \cup \alpha_{ij} \cup C_j$. Let $T_i : S \rightarrow S$ and $T_{jk} = T_{kj} : S \rightarrow S$ denote right Dehn twist maps supported on regular neighborhoods on $S$ of $C_i$ and $C_{jk}$. Let $t_i$ and $t_{jk} = t_{kj}$ be the isotopy classes of $T_i$ and $T_{jk} = T_{kj}$. Then, we have the lantern relation:

$$t_0 \cdot t_1 \cdot t_2 \cdot t_3 = t_{12} \cdot t_{23} \cdot t_{31} = t_{23} \cdot t_{31} \cdot t_{12} = t_{31} \cdot t_{12} \cdot t_{23}.$$

\noindent Since $C_1$, $C_2$, and $C_3$ are boundary components of $S$, it follows that $t_i$ is equal to the identity element of $\Gamma^*(S)$, $i = 1, 2, 3$. Hence: 

$$t_0 = t_{12} \cdot t_{23} \cdot t_{31} = t_{23} \cdot t_{31} \cdot t_{12} = t_{31} \cdot t_{12} \cdot t_{23}.$$

\begin{figure}[!hbp]
\begin{center}
\scalebox{0.50}{\includegraphics{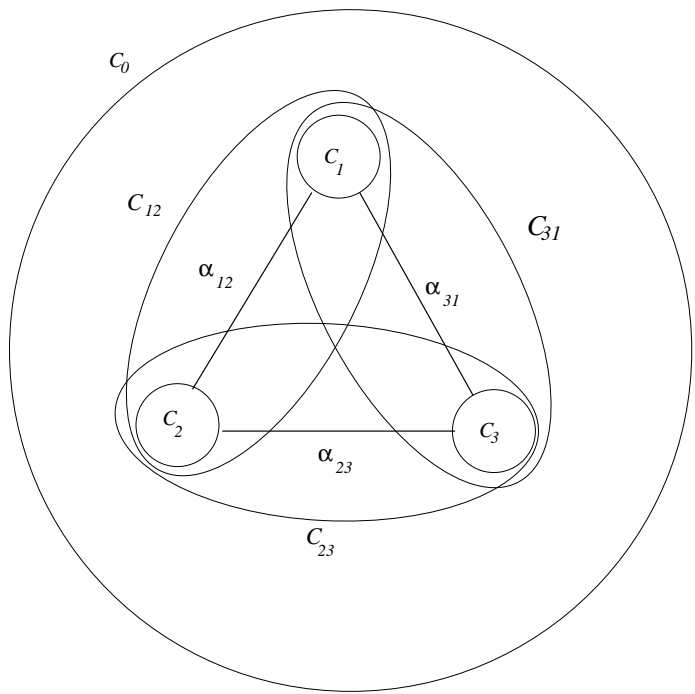}}
\end{center}
\caption{The lantern relation: $t_0 \cdot t_1 \cdot t_2 \cdot t_3 = t_{12} \cdot t_{23} \cdot t_{31}$, where $t_i$ ($t_{jk}$) is the isotopy class 
of a twist map about $C_i$ ($C_{jk}$).}
\label{fig:lantern}\end{figure}

\noindent Suppose that $H : S \rightarrow S$ is orientation-reversing. Let $h \in \Gamma^*(S)$ be the isotopy class of $H : S \rightarrow S$. Since $H : S \rightarrow S$ is orientation-reversing, $ht_{ij}h^{-1}$ is equal to $s_{ij}$ where $s_{ij}$ is a left Dehn twist about $H(C_{ij})$. On the other hand, since $H$ preserves the isotopy class of every essential curve on $S$, $H(C_{ij})$ is isotopic on $S$ to $C_{ij}$. It follows that $s_{ij} = t_{ij}^{-1}$. Likewise, if $C_0$ is essential on $S$, then $h \cdot t_0 \cdot h^{-1} = t_0^{-1}$. On the other hand, if $C_0$ is a boundary component of $S$, then $t_0$ is equal to the identity element of $\Gamma^*(S)$ and, hence, 
$h \cdot  t_0  \cdot h^{-1} = t_0^{-1}$. In any case, $h  \cdot t_0  \cdot  h^{-1} = t_0^{-1}$. We conclude that: 

$$t_0^{-1} = h \cdot t_0 \cdot h^{-1} = h \cdot (t_{12} \cdot t_{23} \cdot t_{31}) \cdot h^{-1} = t_{12}^{-1} \cdot t_{23}^{-1} \cdot t_{31}^{-1}.$$

\noindent This implies that: 

$$t_0 = t_{31} \cdot t_{23} \cdot t_{12}.$$ 

\noindent It follows from the above equations that: 

$$t_{31} \cdot t_{12} \cdot t_{23} = t_{31} \cdot t_{23} \cdot t_{12}$$

\noindent which implies that: 

$$t_{12} \cdot t_{23} = t_{23} \cdot t_{12}.$$

\noindent Hence, the right Dehn twists $t_{12}$ and $t_{23}$ about the essential curves $C_{12}$ and $C_{23}$ on $S$ commute. It follows from Lemma 4.3 of \cite{mcc3} that the geometric intersection $i(C_{12},C_{13}) = 0$. Since this geometric intersection number is equal to two, this is a contradiction. Hence, $H$ is orientation-preserving. 

The case where the genus of $S$ is positive is handled similarly, by using a torus with one hole on $S$ instead of a sphere with four holes on $S$. 
Suppose that $S$ has positive genus. Then there exist transverse essential curves $\alpha$ and $\beta$ on $S$ such that $\alpha$ and $\beta$ have exactly one point of intersection. Let $T_\alpha : S \rightarrow S$ be a homeomorphism representing $t_\alpha$ and $\gamma$ be the image of $\beta$ under $T_\alpha$. Then, since $T_\alpha$ is orientation-preserving:

$$t_\alpha \cdot t_\beta \cdot t_\alpha^{-1} =  t_\gamma.$$

\begin{figure}[!hbp]
\begin{center}
\scalebox{0.50}{\includegraphics{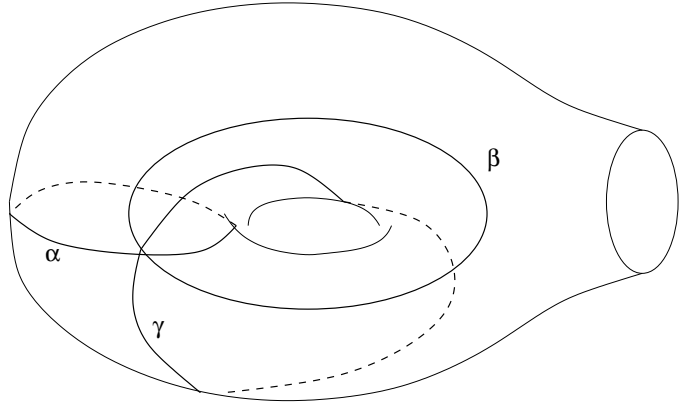}}
\end{center}
\caption{A conjugation relation: $t_\alpha \cdot t_\beta \cdot t_\alpha^{-1} = t_\gamma$}
\label{fig:conjugation}\end{figure}

\noindent Since $\beta$ is essential on $S$ and $T_\alpha : S \rightarrow S$ is a homeomorphism, $\gamma$ is essential on $S$. As before, by conjugating by $h$, we conclude that: 

$$t_\alpha^{-1} \cdot t_\beta^{-1} \cdot t_\alpha =  t_\gamma^{-1}.$$

\noindent This implies that: 

$$t_\alpha^{-1} \cdot t_\beta \cdot t_\alpha =  t_\gamma.$$

\noindent It follows from the above equations that: 

$$t_\alpha \cdot t_\beta \cdot t_\alpha^{-1} =  t_\alpha^{-1} \cdot t_\beta \cdot t_\alpha$$

\noindent which implies: 

$$t_\alpha^2 \cdot t_\beta =  t_\beta \cdot t_\alpha^2.$$

\noindent As before, it follows from Lemma 4.3 of \cite{mcc3} that $i(\alpha,\beta) = 0$. Since $i(\alpha,\beta) = 1$, this is a contradiction. Hence, $H$ is orientation-preserving. 

In any case, $H$ is orientation-preserving.

\end{proof}

A {\it subsurface of $S$} is a surface with boundary $X$ contained in $S$ such that every boundary component of $X$ is either a boundary component of $S$ or disjoint from the boundary of $S$.

\begin{defn} A {\it domain on $S$} is a connected compact subsurface $X$ of $S$ which is not equal to $S$ and each of whose boundary components is either contained in $\partial S$ or is essential on $S$.\label{defn:domain} \end{defn} 

The {\it peripheral boundary components of $X$} are those which are contained in $\partial S$. The remaining boundary components of $X$ are essential. 

\begin{prop} Let $X$ be a domain on $S$. Then $X$ is not a disk; no boundary component of $X$ bounds a disk on $S$; there does not exist an annulus on $S$ whose boundary is equal to the union of a boundary component of $X$ with a boundary component of $S$; and $X$ has at least one essential boundary component on $S$.  
\label{prop:basic1} \end{prop} 

We say that a domain on $S$ is {\it peripheral} if it has at least one peripheral boundary component. We say that a domain on $S$ is {\it monoperipheral} if it has exactly one peripheral boundary component, and {\it biperipheral} if it has exactly two peripheral boundary component. More generally, we say that a domain on $S$ is {\it $k$-peripheral} if it has exactly $k$ peripheral boundary components. 

Let $C$ be a curve on $S$. A {\it regular neighborhood of $C$ on $S$} is an annulus $R$ in the interior of $S$ such that $C$ is a curve on $R$ which does not bound a disc on $R$. Note that a regular neighborhood of a curve on $S$ is a domain on $S$ if and only if the curve is an essential curve on $S$. 

\begin{prop} Let $X$ be a domain on $S$. Let $\alpha$ and $\beta$ be curves on $X$. Then the geometric intersection number of $\alpha$ and $\beta$ in $X$ is equal to the geometric intersection number of $\alpha$ and $\beta$ in $S$.
\label{prop:invgeomint} \end{prop}

\begin{proof} Let $m$ be equal to the geometric intersection number of $\alpha$ and $\beta$ in $X$ and $n$ be equal to the geometric intersection number of $\alpha$ and $\beta$ in $S$. Without loss of generality, we may assume that $\alpha$ and $\beta$ are transverse with exactly $m$ points of intersection. Then there does not exist a disk on $X$ whose boundary is the union of an arc of $\alpha$ and an arc of $\beta$. Since $\alpha$ and $\beta$ meet transversely at $m$ points, $m \geq n$. Suppose that $m > n$. Then there exists a disk $D$ on $S$ whose boundary $\partial_D$ is the union of an arc of $\alpha$ and an arc of $\beta$. 

Let $C$ be a component of $\partial X$. Since $C$ is connected and disjoint from $\partial D$, $C$ is either contained in the interior of $D$ or in the complement of $D$ in $S$. Since $C$ does not bound a disk on $S$, it follows that $C$ is in the complement of $D$. It follows that $D$ is disjoint from $\partial X$. Since $D$ is connected and disjoint from $\partial X$, $D$ is either contained in the interior of $X$ or the complement of $X$ on $S$. Since $\partial D$ is contained in the interior of $X$, it follows that $D$ is contained in the interior of $X$.  Hence, $D$ is a disk on $X$ whose boundary $\partial_D$ is the union of an arc of $\alpha$ and an arc of $\beta$. This is a contradiction. Thus, $m \leq n$ and, hence, $m = n$. 

\end{proof}

The above proof of Proposition \ref{prop:invgeomint} uses a standard  disk argument. Alternatively, one can prove Proposition \ref{prop:invgeomint}, using hyperbolic geometry.

\begin{defn} We say that a domain $X$ on $S$ is {\it elementary} if it is either an annulus or a pair of pants on $S$. 
\label{defn:elemdom} \end{defn}

\begin{prop} Let $X$ be a domain on $S$. Then $X$ is a nonelementary domain on $S$ if and only if there exist curves $\alpha$ and $\beta$ on $S$ such that $i(\alpha,\beta) \neq 0$ and $\alpha$ and $\beta$ are contained in the interior of $X$. \label{prop:nonelemdom} \end{prop}

\begin{proof} Suppose, on the one hand, that $X$ is elementary and $\alpha$ and $\beta$ are curves on $S$ contained in the interior of $X$. Then $\alpha$ is sotopic on $S$ to a boundary components $\partial$ of $X$. Since $\beta$ is in the interior of $X$, $\partial$ and $\beta$ are disjoint. Hence, $i(\alpha,\beta) = i(\partial,\beta) = 0$. 

Suppose, on the other hand, that $X$ is nonelementary. Since $X$ is a domain on $S$, $X$ has at least one essential boundary component on $S$. 

Suppose on the one hand that the genus of $X$ is positive. Then there exists an embedded torus with one hole $Y$ in $X$ which is either equal to $X$ or is a domain on $X$. Let $\alpha$ and $\beta$ be curves on $Y$ which intersect transversally and have exactly one point of intersection. Then, by Proposition \ref{prop:invgeomint},  $i_S(\alpha,\beta) = i_X(\alpha,\beta) = i_Y(\alpha,\beta) = 1$.  

Suppose that the genus of $X$ is zero. Since $X$ is a domain on $S$, $X$ has at least one essential boundary component on $S$. In particular, $X$ is not a disc. Since $X$ is nonelementary, $X$ is not an annulus or a pair of pants. Hence, $X$ has at least four boundary components. Thus there exists an embedded sphere with four holes $Y$ in $X$ which is either equal to $X$ or is a domain on $X$. Let $\alpha$ and $\beta$ be curves on $Y$ which intersect transversally, have exactly two points of intersection, and are such that the complement of $\alpha \cup \beta$ in $Y$ has exactly four 
components, each of which contains exactly one boundary component of $Y$. Then, by Proposition \ref{prop:invgeomint},  $i_S(\alpha,\beta) = i_X(\alpha,\beta) = i_Y(\alpha,\beta) = 2$.
 
\end{proof}

\begin{prop} Let $X$ be a domain on $S$ and $Y$ be a subsurface of $X$.  If $Y$ is a domain on $X$, then $Y$ is a domain on $S$. 
\label{prop:transitivity} \end{prop}

\begin{proof} By Proposition \ref{prop:basic1}, $S$ is not a disk and no boundary component of $X$ bounds a disk on $S$. 

Clearly $Y$ is a compact, connected, orientable subsurface of $S$ which is not equal to $S$. Let $\partial$ be a boundary component of $Y$. Since $Y$ is a domain on $X$, $\partial$ is either a boundary component of $X$ or an essential curve on $X$. 

Suppose, on the one hand, that $\partial$ is a boundary component of $X$. Since $X$ is a domain on $S$, it follows that $\partial$ is either a boundary component of $S$ or an essential curve on $S$. 

Suppose, on the other hand, that $\partial$ is an essential curve on $X$. 
In particular, $\partial$ is in the interior of $X$. Hence, $\partial$ is in the interior of $S$. 

Suppose that $\partial$ is not an essential curve on $S$. Then $\partial$ either bounds a disk $D$ on $S$ or cobounds an annulus $A$ on $S$ with a boundary component $\epsilon$ of $S$.  

Suppose that $\partial$ bounds a disk $D$ on $S$. Since no boundary component of $X$ bounds a disk on $S$, no boundary component of $X$ can be contained in $D$. Hence, $D$ is disjoint from the boundary of $X$ but intersects the interior of $X$, since it contains $\partial$. Since $D$ is connected, this implies that $D$ is contained in $X$. Since $\partial$ is an essential curve on $X$, this is a contradiction. Hence, $\partial$ does not bound a disk on $S$. 

Hence, $\partial$ cobounds an annulus $A$ on $S$ with a boundary component $\epsilon$ of $S$. Since no boundary component of $X$ can bound a disk on $S$ or cobound an annulus with the boundary component $\epsilon$ of $S$, no boundary component of $X$ can be contained in the interior of $A$. Hence, 
$A \setminus \epsilon$ is disjoint from the boundary of $X$ but intersects the interior of $X$, since it contains $\partial$. Since $A \setminus \epsilon$ is connected and $X$ is compact, this implies that $A$ is contained in $X$. As before, this is a contradiction, and, hence, $\partial$ is an essential curve on $S$.

This shows that each boundary component of $Y$ is either a boundary component of $S$ or an essential curve on $S$. That is to say, $Y$ is a domain on $S$. 

\end{proof}

\begin{prop} Let $X$ be a domain on $S$ and $\alpha$ be an essential curve on $X$. Then there exists an essential curve $\beta$ on $X$ such that the geometric intersection number of $\alpha$ and $\beta$ on $S$ is not equal to zero.\label{prop:geomint} \end{prop} 

\begin{proof} Since $X$ is a domain on $S$, $X$ is a connected, compact, orientable surface with boundary. Since $\alpha$ is an essential curve on $X$, it follow from Proposition \ref{prop:detectisotopy}, that there exists an essential curve $\beta$ on $X$ such that the geometric intersection number $i_X(\alpha,\beta)$ of $\alpha$ and $\beta$ on $X$ is not equal to zero. Since $X$ is a domain on $S$ and $\alpha$ and $\beta$ are curves on $X$, it follows from Proposition \ref{prop:invgeomint} that $i_X(\alpha,\beta)$ is equal to the geometric intersection number $i_S(\alpha,\beta)$ of $\alpha$ and $\beta$ on $S$. This implies that $i_S(\alpha,\beta) = i_X(\alpha,\beta) \neq 0$, completing the proof.
\end{proof}

\begin{prop} Let $X$ and $Y$ be domains on $S$. Suppose that $X$ is isotopic on $S$ to a domain on $Y$. Then $Y$ is not isotopic on $S$ to a domain on $X$.
\label{prop:partialorder} \end{prop}

\begin{proof} Suppose that $Y$ is isotopic on $S$ to a domain on $X$. Let $X_1$ be a domain on $Y$ such that $X$ is isotopic to $X_1$ on $S$. Let $Y_1$ be a domain on $X$ such that $Y$ is isotopic to $Y_1$ on $S$. Since $Y$ is isotopic to $Y_1$ on $S$ and $Y_1$ is a domain on $X$, it follows that $X$ is isotopic to a domain $X_2$ on $S$ such that $Y$ is a domain on $X_2$. 

Since $X_1$ is a domain on $Y$ and $Y$ is a domain on $X_2$, it follows from Proposition \ref{prop:transitivity} that $X_1$ is a domain on $X_2$. 

Since $X_1$ is a domain on $X_2$, $X_1$ has an essential boundary component $\alpha$ on $X_2$. Since $\alpha$ is an essential curve on $X_2$, it follows from Proposition \ref{prop:geomint} that there exists an essential curve $\beta$ on $X_2$ such that the geometric intersection number of $\alpha$ and $\beta$ on $S$ is not equal to zero.

Since $X_1$ is isotopic to $X_2$ on $S$, the essential boundary component $\alpha$ of $X_1$ on $S$ is isotopic to an essential boundary component of $X_2$ on $S$. Hence, $\alpha$ is isotopic to a curve $\alpha_1$ on $S$ such that $\alpha_1$ is disjoint from $X_2$. Since $\alpha$ is isotopic to $\alpha_1$ on $S$, it follows that $i(\alpha,\beta) = i(\alpha_1,\beta)$.
Since $\alpha_1$ is disjoint from $X_2$ and $\beta$ is contained in $X_2$, $\alpha_1$ is disjoint from $\beta$ and, hence, $i(\alpha_1,\beta) = 0$. Hence, $i(\alpha,\beta) = 0$, which is a contradiction. 

Hence, $Y$ is not isotopic on $S$ to a domain on $X$.

\end{proof} 

The following is a weak converse for Proposition \ref{prop:transitivity}.

\begin{prop} Let $X$ be a domain on $S$ and $Y$ be a subsurface of $X$.  If 
$Y$ is a domain on $S$, then $Y$ is isotopic on $S$ to either $X$, or a domain on $X$, or a regular neighborhood of an essential boundary component of $X$. 
\label{prop:weakconv} \end{prop}

\begin{proof} Assume that $Y$ is not isotopic on $S$ to either $X$ or a regular neighborhood of an essential boundary component of $X$. 

Let $\partial$ be a boundary component of $Y$. We may assume that $\partial$ is not a boundary component of $X$. Hence, $\partial$ is not a boundary component of $S$. Since $Y$ is a domain on $S$, it follows that $\partial$ is an essential curve on $S$.

Suppose that $\partial$ is not an essential curve on $X$. 

Since $\partial$ is an essential curve on $S$ it cannot bound a disk on $S$. Hence, it cannot bound a disk on $X$. 

Hence, $\partial$ must cobound an annulus $A$ on $X$ with a boundary component $\epsilon$ of $X$. 

Since $\partial$ is essential on $S$ and $A$ is an annulus on $S$, $\epsilon$ is not a boundary component of $S$. 

Hence, since $X$ is a domain on $S$, $\epsilon$ is an essential boundary component of $X$ on $S$. 

Since $Y$ is contained in $X$ and the interior of $Y$ is disjoint from the complement of  $\partial \cup \epsilon$ in $A$, the interior of $Y$ is contained in either the interior of $A$ or the complement of $A$ in $X$. Hence, $Y$ is contained in either 
$A$ or the closure of the complement of $A$ in $X$.

Suppose that $Y$ is contained in $A$. Since $Y$ is a domain on $S$ and $A$ is contained in the interior of $S$, it follows that $Y$ is isotopic on $S$ to $A$. Since $Y$ is not isotopic to a regular neighborhood of an essential boundary component of $X$ on $S$, this is a contradiction.

Hence, $Y$ is contained in the closure of the complement of $A$ in $X$. 

Let $Y' = Y \cup A$. Note that $Y'$ is a domain on $S$ which is isotopic to $Y$ on $S$, is contained in $X$, and has one less boundary component in the interior of $X$ than $Y$. 

It follows by induction, that $Y$ is isotopic to a domain on $S$ which is contained in $X$ and which has all of its boundary components in the boundary of $X$. In other words, $Y$ is isotopic on $S$ to $X$. 

\end{proof}

\begin{defn} Let $\mathcal{F}$ be a collection of pairwise disjoint domains on $S$. Let $Y$ be the closure of the complement of $|\mathcal{F}|$ in $S$. The {\it codomains of $\mathcal{F}$} are the components of $Y$. 
\label{defn:codomains} \end{defn}

Note that the codomains of a collection of pairwise disjoint domains on $S$ are themselves domains on $S$.  

\begin{defn} A nonempty collection of pairwise disjoint domains on $S$ is called a {\it system of domains on $S$} if the domains in the collection are pairwise nonisotopic.\label{defn:system} \end{defn} 

The collection of codomains of a system of domains on $S$ is a collection of pairwise disjoint domains on $S$. However, the collection of codomains of a system of domains on $S$ is not necessarily a system of domains on $S$, since two distinct codomains of a system of domains on $S$ may be isotopic. 

\section{Preliminaries on simplicial complexes} 

In this section, we discuss some terminology from the theory of abstract simplicial complexes. A reference for this material is Munkres \cite{munkres}.

\begin{defn} Let $V$ be a set. An {\it abstract simplicial complex $K$ with vertex set $V$} is a collection of finite subsets of $V$ such that: 
\begin{enumerate} 
\item \label{abstract1} if $v \in V$, then $\{v\} \in K$; 
\item \label{abstract2} if $\sigma$ is an element of $K$ and $\tau$ is a subset of $\sigma$, then $\tau$ is an element of $K$. 
\end{enumerate}
\label{defn:abstract} \end{defn}

\noindent The term {\it simplicial complex} in this paper shall refer to an abstract simplicial complex. 

Let $K$ be a simplicial complex with vertex set $V$. If $x$ is an element of $V$, then we say that $x$ is a {\it vertex of $K$}. If $\sigma$ is an element of $K$ and $x$ is an element of $\sigma$, then we say that $\sigma$ is a {\it simplex of $K$} and $x$ is a {\it vertex of $\sigma$}. Note that each vertex of each simplex of $K$ is a vertex of $K$. If $\sigma$ has $k + 1$ vertices, then we say that $\sigma$ is a {\it $k$-simplex of $K$}. If $x$ is an element of $V$, then we also say that the corresponding $0$-simplex $\{x\}$ of $K$ is a vertex of $K$. If $e$ is a $1$-simplex of $K$, then we say that {\it $e$ is an edge of $K$}. If $\Delta$ is a $2$-simplex of $K$, then we say that {\it $\Delta$ is a triangle of $K$}. 

\begin{defn} Let $K$ be an abstract simplicial complex. Let $F$ be a subcollection of $K$. We say that {\it $F$ is a subcomplex of $K$} if each subset $\tau$ of an element of $F$ is an element of $F$.
\label{defn:subcomplex} \end{defn}

\noindent Let $F$ be a subcomplex of an abstract simplicial complex $K$. Note that $F$ is itself an abstract simplicial complex, and the vertex set of $F$ is a subset of the vertex set of $K$. 
 
\begin{prop} Let $K$ be a simplicial complex with vertex set $V$ and $W$ be a subset of $V$. Let $K_W$ be the set of all simplices of $K$ which have all of their vertices in $W$. Then $K_W$ is a subcomplex of $K$ with vertex set $W$. \label{prop:inducedsubcomplex} \end{prop}

\begin{defn} Let $K$, $W$, and $K_W$ be as in Proposition \ref{prop:inducedsubcomplex}. We say that {\it $K_W$ is the subcomplex of $K$ induced by the subset $W$ of the set of vertices $V$ of $K$}.\label{defn:inducedsubcomplex} \end{defn}

\noindent Note that the subcomplex of a simplicial complex induced by a subset of its vertices is itself a simplicial complex. Moreover, it is completely determined by the simplicial complex $K$ and its vertex set $W$.

Let $K$ be a simplicial complex. For each nonnegative integer $n$, the {\it $n$-skeleton $K_n$ of $K$} is the subcomplex of $K$ consisting of all $k$-simplices of $K$ with $k \leq n$. Note that $V$ is equal to the support $|K_n|$ of $K_n$ (i.e. the union of all the $k$-simplices of $K$ with 
$k \leq n$).

Note that if $F$ is a subcomplex of an abstract simplicial complex $K$ and $n$ is a nonnegative integer, then the $n$-skeleton $F_n$ of $F$ is a subcomplex of the $n$-skeleton $K_n$ of $K$. 

If $\tau \subset \sigma \in K$, then we say that {\it $\tau$ is a face of $\sigma$}. A {\it maximal simplex of a simplicial complex $K$} is a simplex which is not a proper face of any simplex of $K$. $K$ is {\it finite-dimensional} if there exists an integer $N$ such that every simplex of $K$ is a $k$-simplex for some $k \leq N$. If $K$ is 
finite-dimensional, then the {\it dimension of  $K$} is the minimum such integer $N$. If the dimension of $K$ is $N$, then a {\it top-dimensional simplex of $K$} is an $N$-simplex of $K$.

\begin{defn} A simplicial complex $K$ is a {\it flag complex} if the following holds: 
\begin{itemize} 
\item If $\{x_0,\ldots,x_n\} \subset K_0$ such that $\{x_i,x_j\}$ is an edge of $K$ for 
$0 \leq i < j \leq n$, then $\{x_0,\ldots,x_n\}$ is a simplex of $K$. 
\end{itemize} 
\label{defn:flag} \end{defn}

\begin{defn} Let $\alpha$, $\beta$, and $\delta$ be simplices of a simplicial complex $K$. We say that $\alpha$ is joined to $\beta$ by $\delta$ if $\delta = \alpha \cup \beta$. \label{defn:joined} \end{defn}

Note that if $\alpha$ is joined to $\beta$ by simplices $\delta$ and $\epsilon$ of $K$, then $\delta = \epsilon$.

Note also that two simplices of a simplicial complex are joined by a vertex of that simplicial complex if and only if they are both equal to that vertex.  

\begin{defn} Let $\alpha$ and $\beta$ be simplices of a simplicial complex $K$ which are joined in $K$ by a simplex $\delta$ of $K$. Then we say that {\it $\delta$ is the join of $\alpha$ and $\beta$ in $K$}.\label{defn:join} \end{defn}

\begin{defn}[the star of a vertex of a simplicial complex] Let $x$ be a vertex of a simplicial complex $K$. The {\it star of $x$ in $K$} is the subcomplex $St(x,K)$ of $K$ whose simplices are the simplices of $K$ which contain the vertex $x$ together with all the faces of such simplices of $K$.\label{defn:star} \end{defn}

Let $K$ be a simplicial complex and $x$ be a vertex of $K$. Note that the $0$-skeleton $St_0(x,K)$ of $St(x,K)$ is the set of all vertices $w$ of $K$ such that $\{x,w\}$ is a simplex of $K$. 

\begin{prop} Let $K$ be a flag complex. Let $x$ and $y$ be vertices of $K$. Then the following are equivalent: 

\begin{enumerate} 

\item \label{flagstar1}  $St(x,K) = St(y,K)$.

\item \label{flagstar2}  $St_0(x,K) = St_0(y,K)$.

\end{enumerate} 

\label{prop:flagstar} \end{prop}

\begin{proof} Clearly (\ref{flagstar1}) implies (\ref{flagstar2}). 

We shall now show that (\ref{flagstar2}) implies (\ref{flagstar1}). To this end, let $\tau$ be a simplex of $St(x,K)$. We need to show that $\tau$ is a simplex of $St(y,K)$. In other words, we need to show that $\tau \cup \{y\}$ is a simplex of $K$.

Since $K$ is a flag complex, it suffices to show that any two distinct vertices of $\tau \cup \{y\}$ are joined by an edge of $K$. To this end, let $w$ and $z$ be distinct vertices of $\tau \cup \{y\}$. If $w$ and $z$ are both vertices of the simplex $\tau$ of $K$, then $\{w,z\}$ is an edge of the simplex $\tau$ of $K$ and, hence, an edge of $K$. 

Hence, we may assume that $w$ is a vertex of $\tau$ and $z = y$. 
Since $w$ is a vertex of $\tau$ and $\tau$ is a simplex of $St(x,K)$, $w$ is a vertex of $St(x,K)$ and, hence, of $St(y,K)$. Since $w$ is a vertex of $St(y,K)$, $\{w,y\}$ is a simplex of $K$. Since $z = y$, we conclude that $\{w,z\}$ is a simplex of $K$. 

Hence, every two distinct vertices of $\tau \cup \{y\}$ are joined by an edge of $K$. Since $K$ is a flag complex, this implies that $\tau \cup \{y\}$ is a simplex of $K$. In other words, $\tau$ is a simplex of $St(y,K)$. 

This proves that $St(x,K)$ is a subcomplex of $St(y,K)$. By a symmetric argument, it follows that $St(y,K)$ is a subcomplex of $St(x,K)$. 

This proves that $St(x,K) = St(y,K)$.

\end{proof}

\begin{defn}[the link of a vertex of a simplicial complex] Let $x$ be a vertex of a simplicial complex $K$. The  {\it link of $x$ in $K$} is the subcomplex $Lk(x,K)$ of $K$ whose simplices are the simplices of $St(x,K)$ which do not have $x$ as a vertex. \label{defn:linksimplex} \end{defn}

Suppose that $x$ is a vertex of a simplicial complex $K$. Note that $\{x\}$ is joined to each of the simplices of $Lk(x,K)$; and the simplices of $St(x,K)$ are precisely the faces of the joins of $\{x\}$ with the simplices of $Lk(x,K)$.  

\begin{defn} Let $K$ be a simplicial complex with vertex set $V$ and $L$ be a simplicial complex with vertex set $W$. A {\it simplicial map from $K$ to $L$} is a map $\varphi : V \rightarrow W$ such that, for each simplex $\sigma$ of $K$, $\varphi(\sigma)$ is a simplex of $L$. \label{defn:morphism} \end{defn}

\noindent Let $\varphi : V \rightarrow W$ be a simplicial map from a simplicial complex $K$ with vertex set $V$ to a simplicial complex $L$ with vertex set $W$. Since $\varphi : V \rightarrow W$ is a simplicial map from $K$ to $L$, the rule $\sigma \mapsto \varphi(\sigma)$ determines a map from $K$ to $L$. We denote this map by $\varphi : K \rightarrow L$ and we say that {\it $\varphi : K \rightarrow L$ is a simplicial map from $K$ to $L$}. In order to distinguish $\varphi : V \rightarrow W$ and $\varphi : K \rightarrow L$, we shall say that $\varphi : V \rightarrow W$ is the {\it vertex correspondence associated to $\varphi : K \rightarrow L$}.

Note that the map $\varphi : K \rightarrow L$ is both determined by and determines the map $\varphi : V \rightarrow W$. The map $\varphi : K \rightarrow L$ is injective if and only if $\varphi : V \rightarrow W$ is injective. If $\varphi: K \rightarrow L$ is surjective, then $\varphi : V \rightarrow W$ is surjective. The converse, however, is not necessarily true. For instance, if $L$ has at least one edge $e$ and $K$ is equal to the zero skeleton $L_0$ of $L$, then the vertex set $V$ of $K$ is equal to the vertex set $W$ of $L$, the identity map $\varphi : V \rightarrow W$ is surjective, but the corresponding map $\varphi: K \rightarrow L$ is not surjective, since the edge $e$ of $L$ is not in the image of $\varphi : K \rightarrow L$.

If $\varphi : K \rightarrow L$ is a simplicial map, then $\varphi(K)$ is a subcomplex of $L$. 

\begin{defn} Let $K$ and $L$ be abstract simplicial complexes. A simplicial isomorphism $\varphi : K \rightarrow L$ is a simplicial map $\varphi: K \rightarrow L$ for which there exists a simplicial map $\psi: L \rightarrow K$ such that $\varphi : K \rightarrow L$ and $\psi : L \rightarrow K$ are inverse functions. In the case where $K = L$, we call a simplicial isomorphism $\varphi: K \rightarrow L$ a simplicial automorphism of $K$. \label{defn:isomorphism} \end{defn} 

Note that a simplicial map $\varphi: K \rightarrow L$ is a simplicial isomorphism if and only if $\varphi: K \rightarrow L$ is bijective. By the previous observations, if $\varphi: K \rightarrow L$ is bijective, then $\varphi : V \rightarrow W$ is bijective. The converse however need not be true, as illustrated by the inclusion $\varphi : L_0 \rightarrow L$ of the zero skeleton $L_0$ of a complex $L$ with at least one edge $e$, whose corresponding vertex correspondence $\varphi : W \rightarrow W$ is the identity map of the set of vertices $W$ of $L$.

\begin{prop} Let $K$ be a simplicial complex, $\varphi \in Aut(K)$, and $x$ be a vertex of $K$. Then $\varphi(St(x,K)) = St(\varphi(x),K)$ and $\varphi(Lk(x,K)) = Lk(\varphi(x),K)$.
\label{prop:invariance} \end{prop}

\section{Exchange automorphisms of simplicial complexes} 

We shall need the following basic results about automorphisms of abstract simplicial complexes. 

\begin{defn} Let $K$ be a simplicial complex, $\{x,y\}$ be a pair of vertices of $K$, and $\varphi : K \rightarrow K$ be an automorphism of $K$. We say that {\it $\varphi$ is a simple exchange of $K$ exchanging the vertices $x$ and $y$ of $K$} if $\varphi(x) = y$, $\varphi(y) = x$, and $\varphi(z) = z$ for every vertex $z$ of $K$ which is neither equal to $x$ nor equal to $y$.\label{defn:simpleexchange} \end{defn} 

Let $\varphi : K \rightarrow K$ be a simple exchange of a simplicial complex $K$ exchanging the vertices $x$ and $y$ of $K$. Note that $\varphi : K \rightarrow K$ is equal to the identity map $id_K : K \rightarrow K$ of $K$ if and only if $x = y$. In this case, we say that {\it $\varphi$ is a trivial simple exchange}. Let $K$ be a simplicial complex, $\varphi : K \rightarrow K$ be a simple exchange of $K$ exchanging the vertices $x$ and $y$ of $K$, and $\psi : K \rightarrow K$ be a simple exchange of $K$ exchanging the vertices $u$ and $v$ of $K$. Then $\varphi = \psi$ if and only if either $x = y$ and $u = v$ or $\{x,y\} = \{u,v\}$. In particular, a nontrivial simple exchange of $K$ exchanges a unique pair of distinct vertices of $K$. 

\begin{exam} Let $K(n)$ denote the simplicial complex of all subsets of the set $\{1,\ldots,n\}$. Then, for every pair of distinct vertices, $i$ and $j$, of $K(n)$, the standard transposition $(i,j)$ in the group of permutations $\Sigma_n$ of $\{1,\ldots,n\}$ extends to a simple exchange of $K(n)$ which exchanges $i$ and $j$. These simple exchanges, of course, generate the group of simplicial automorphisms of $K(n)$, which is naturally isomorphic to the symmetric group $\Sigma_n$, the group of permutations of the vertex set $\{1,\ldots,n\}$ of $K(n)$. 
\label{exam:symmetricgroup} \end{exam}

\begin{defn} Let $K$ be a simplicial complex. Let $x$ and $y$ be distinct vertices of $K$. We say that {\it $x$ and $y$ are exchangeable in $K$} if there exists a simple exchange of $K$ exchanging $x$ and $y$.\label{defn:exchangable} \end{defn}

\begin{figure}[!hbp]
\begin{center}
\scalebox{0.60}{\includegraphics{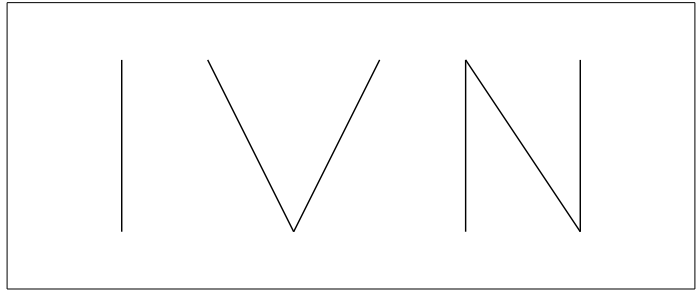}}
\end{center}
\caption{Three complexes: one with all vertices exchangeable; one with some vertices but not all exchangeable; and one with no vertices exchangeable.}
\label{fig:ivn}\end{figure}
  
\begin{figure}[!hbp]
\begin{center}
\scalebox{0.60}{\includegraphics{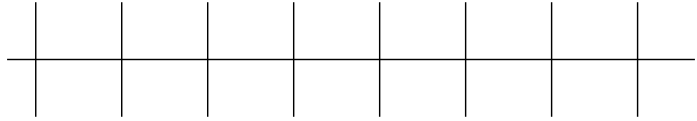}}
\end{center}
\caption{Example \ref{exam:lineofedges}: a line of edges.}
\label{fig:lineofedges}\end{figure}

\begin{exam} Let $V = \mathbb{Z} \times \{-1,0,1\}$ and $K$ be the 
one-dimensional simplicial complex on $V$, illustrated in Figure \ref{fig:lineofedges}, whose edges are the pairs $\{(m,0),(m + 1,0)\}$ and $\{(m,0),(m,\epsilon)\}$ with $m \in \mathbb{Z}$ and $\epsilon \in \{-1,1\}$. Then two distinct vertices $x$ and $y$ of $K$ are exchangeable if and only if $\{x,y\} = \{(m,-1),(m,1)\}$ for some $m \in \mathbb{Z}$. 

Note that for any subset $W$ of $\mathbb{Z}$ there is a unique automorphism $\varphi_W : K \rightarrow K$ such that $\varphi_W(m,t)$ is equal to $(m,-t)$ if $m \in W$ and $(m,t)$ otherwise. In particular, $\varphi_\emptyset = id_K : K \rightarrow K$. If $U$ and $V$ are subsets of $\mathbb{Z}$, then $\varphi_U \circ \varphi_V = \varphi_{U \triangle V}$. In particular, $\varphi_W \circ \varphi_W = id_K : K \rightarrow K$. It follows that the collection $\{ \varphi_W | W \subset \mathbb{Z}\}$ of automorphisms of $K$ is a subgroup $B_K$ of the group of automorphisms of $K$, $Aut(K)$, naturally isomorphic to the Boolean algebra $\mathcal{B}(\mathbb{Z})$ of all subsets of $\mathbb{Z}$.  

Let $D_K$ be the group of automorphisms of $K$ generated by the translation $(m,n) \mapsto (m + 1, n)$ and the involution  $(m,n) \mapsto (1-m,n)$. Note that this involution has no fixed vertices in $K$. The subgroup $D_K$ of $Aut(K)$ is naturally isomorphic to the infinite dihedral group $D_\infty$ of isometries of $\mathbb{Z}$ equipped with its standard metric.

The group of automorphism of $K$, $Aut(K)$, is a split extension of its subgroup $D_K$ by its normal subgroup $B_K$, and we have the following commutative diagram.  

\[ \begin{array}{cccccccccc}
1 &\longrightarrow & \mathcal{B}(\mathbb{Z})   & \longrightarrow &
\mathcal{B}(\mathbb{Z}) \rtimes D_{\infty} & \longrightarrow  & D_{\infty}=Isom(\mathbb{Z}) & \longrightarrow & 1\\
\   & \  & \simeq\big\downarrow & \  & \simeq\big\downarrow & \  & \simeq\big\downarrow 
& \ & \  \\
1 & \longrightarrow & B_K &\longrightarrow & Aut(K)  &\longrightarrow &
D_K & \longrightarrow & 1 
\label{diagram1} \end{array}\]

\label{exam:lineofedges} \end{exam}

The following result gives a basic necessary and sufficient condition for two vertices of a simplicial complex to be exchangeable. 

\begin{prop} Let $K$ be a simplicial complex and $x$ and $y$ be vertices of $K$. Let $F$ be the subcomplex of $K$ consisting of all simplices of $K$ which have neither $x$ nor $y$ as a vertex. Then the following are equivalent: 

\begin{enumerate} 

\item \label{exchangecondition1} $x$ and $y$ are exchangeable in $K$.

\item \label{exchangecondition2} $St(x,K) \cap F = St(y,K) \cap F$.

\end{enumerate} 

\label{prop:exchangecondition} \end{prop}

\begin{proof} First, we prove that (\ref{exchangecondition1}) implies (\ref{exchangecondition2}). To this end, suppose that $x$ and $y$ are exchangeable in $K$. Since $x$ and $y$ are exchangeable in $K$, there is a unique automorphism $\varphi : K \rightarrow K$ such that $\varphi(x) = y$, $\varphi(y) = x$, and $\varphi(z) = z$ for every vertex $z$ of $F$.  

Suppose that $\sigma$ is a simplex of $St(x,K) \cap F$. In other words, suppose that $\{x\} \cup \sigma$ is a simplex of $K$ and $\sigma$ is a simplex of $F$. Then $\varphi(\{x\} \cup \sigma) = \{\varphi(x)\} \cup \varphi(\sigma) = \{y\} \cup \sigma$ is a simplex of $K$. Since $\{y\} \cup \sigma$ is a simplex of $K$ and $\sigma$ is a simplex of $F$, it follows that $\sigma$ is a simplex of $St(y,K) \cap F$. This proves that $St(x,K) \cap F \subset St(y,K) \cap F$. Likewise, $St(y,K) \cap F) \subset St(x,K) \cap F$ and, hence, $St(x,K) \cap F = St(y,K) \cap F$. 
This proves that (\ref{exchangecondition1}) implies (\ref{exchangecondition2}).

Now we prove that (\ref{exchangecondition2}) implies (\ref{exchangecondition1}). To this end, suppose that $St(x,K) \cap F = St(y,K) \cap F$. 

Consider the bijection $\varphi : K_0 \rightarrow K_0$ defined by the rule $\varphi(x) = y$, $\varphi(y) = x$, and $\varphi(z) = z$ for every vertex $z$ of $F$. Since $\varphi : K_0 \rightarrow K_0$ is an involution of $K_0$, it suffices to prove that $\varphi$ extends to a simplicial map $\varphi : K \rightarrow K$. In other words, it suffices to prove that $\varphi(\tau)$ is a simplex of $K$ for every simplex $\tau$ of $K$. 

To this end, suppose that $\tau$ is a simplex of $K$. If $x$ and $y$ are both vertices of $\tau$, then $\varphi(\tau)$ is equal to the simplex $\tau$ of $K$. Likewise, if neither $x$ nor $y$ is a vertex of $\tau$, then $\varphi(\tau)$ is equal to the simplex $\tau$ of $K$. 

Suppose that $x$ is a vertex of $\tau$ and $y$ is not a vertex of $\tau$. Let $\sigma = \tau \setminus \{x\}$. Since $x$ is a vertex of $\tau$ and $\sigma = \tau \setminus \{x\}$, it follows that $\tau = \{x\} \cup \sigma$. Since $\tau$ is a simplex of $K$, this implies that $\sigma$ is a simplex of $St(x,K)$. Since $y$ is not a vertex of $\tau$ and $\sigma = \tau \setminus \{x\}$, $\sigma$ is a simplex of $F$. 
This implies that $\sigma$ is a simplex of $St(x,K) \cap F$ and, hence, of $St(y,K) \cap F$. Since $\sigma$ is a simplex of $St(y,K)$, $\{y\} \cup \sigma$ is a simplex of $K$. 

Since $\sigma$ is a simplex of $F$ and $\tau = \{x\} \cup \sigma$, it follows that $\varphi(\tau) = \{\varphi(x)\} \cup \varphi(\sigma) = \{y\} \cup \sigma$. Hence, $\varphi(\tau)$ is a simplex of $K$. 

This shows that if $x$ is a vertex of $\tau$ and $y$ is not a vertex of $\tau$, then $\varphi(\tau)$ is a simplex of $K$. Likewise, if $y$ is a vertex of $\tau$ and $x$ is not a vertex of $\tau$, then $\varphi(\tau)$ is a simplex of $K$. 

In any case, $\varphi(\tau)$ is a simplex of $K$. 

This proves that $\varphi : K_0 \rightarrow K_0$ extends to a simplicial map $\varphi : K \rightarrow K$. Since $\varphi : K_0 \rightarrow K_0$ is an involution, its extension $\varphi : K \rightarrow K$ is an involution. Hence, this extension $\varphi : K \rightarrow K$ is a simplicial automorphism exchanging $x$ and $y$. 

This proves that (\ref{exchangecondition2}) implies (\ref{exchangecondition1}).

\end{proof}

\begin{rem} Note that $Star(x,K)$ joins $x$ to the subcomplex $F \cap Lk(x,K)$ of $F$ and, in the case where $\{x,y\}$ is an edge of $K$, also to $y$. Likewise, $Star(y,K)$ joins $y$ to the subcomplex $F \cap Lk(y,K)$ of $F$ and, in the case where $\{x,y\}$ is an edge of $K$, also to $x$. Roughly speaking, the above {\it exchangeability condition}, condition (\ref{exchangecondition2}), states that $F$ is a sort of hyperplane of reflection across which the vertices $x$ and $y$ of $K$ are able to be reflected since they have been symmetrically joined to $F$ along a subcomplex $G$ of $F$ (i.e. along $F \cap Lk(x,K) = F \cap Lk(y,K)$) and, in the case where $\{x,y\}$ is an edge of $K$, to one another. 
\label{rem:exchangecondition} \end{rem}

The following propositions are refinements of Proposition \ref{prop:exchangecondition} corresponding to the situations where $\{x,y\}$ is or is not an edge of $K$.

\begin{prop} Let $K$ be a simplicial complex and $x$ and $y$ be distinct vertices of $K$ which are not connected by an edge of $K$. Then the following are equivalent: 

\begin{enumerate} 

\item \label{xynotanedge1} $x$ and $y$ are exchangeable in $K$.

\item \label{xynotanedge2} $Lk(x,K) = Lk(y,K)$.

\end{enumerate} 

\label{prop:xynotanedge} \end{prop}

\begin{proof} Let $F$ be the subcomplex of $K$ consisting of all simplices of $K$ which have neither $x$ nor $y$ as a vertex. Since $x$ and $y$ are not joined by an edge of $K$, it follows that $Lk(x,K) = St(x,K) \cap F$ and $Lk(y,K) = St(y,K) \cap F$. 
\end{proof}

\begin{prop} Let $K$ be a simplicial complex and $x$ and $y$ be vertices of $K$ which are connected by an edge of $K$. Let $F$ be the subcomplex of $K$ consisting of all simplices of $K$ which have neither $x$ nor $y$ as a vertex. Suppose that $K$ is a flag complex. Then the following are equivalent: 

\begin{enumerate} 

\item \label{xyanedge1} $x$ and $y$ are exchangeable in $K$.

\item \label{xyanedge2} $St(x,K) = St(y,K)$.

\end{enumerate} 

\label{prop:xyanedge} \end{prop}

\begin{proof} Suppose that $St(x,K) = St(y,K)$. Then $St(x,K) \cap F = St(y,K) \cap F$. It follows from Proposition \ref{prop:exchangecondition} that $x$ and $y$ are exchangeable. This proves that (\ref{xyanedge2}) implies (\ref{xyanedge1}). 

We shall now show that (\ref{xyanedge1}) implies (\ref{xyanedge2}). Suppose that $x$ and $y$ are exchangeable in $K$. It follows from Proposition \ref{prop:exchangecondition} that $St(x,K) \cap F = St(y,K) \cap F$. We must show that $St(x,K) = St(y,K)$. To this end, suppose that $\tau$ is a simplex of $St(x,K)$. Let $\sigma = \tau \cup \{x\}$. Since $\tau$ is a simplex of $St(x,K)$, it follows that $\sigma$ is a simplex of $K$. 

Let $\rho = \sigma \cup \{y\}$. We shall show that $\rho$ is a simplex of $K$. Since $K$ is a flag complex, it suffices to show that any two distinct vertices of $\rho$ are joined by an edge of $K$. To this end, let $w$ and $z$ be vertices of $\rho$. 
If neither $w$ nor $z$ is equal to $y$, then $w$ and $z$ are vertices of the simplex $\sigma$ of $K$ and, hence, are joined by an edge of $K$. Hence, we may assume that $z = y$. This implies that $w$ is not equal to $y$. It follows that $x$ and $w$ are both vertices of the simplex $\sigma$ of $K$. Hence, $w$ is a vertex of $St(x,K)$.
If $w = x$, then $w$ and $z$ are vertices of the simplex $\{x,y\}$ of $K$. 

Hence, we may assume that $w$ is not equal to $x$.

Since $w$ is a vertex of $St(x,K)$ and $w$ is not equal to $x$ or $y$, it follows that $w$ is a vertex of $St(x,K) \cap F$ and, hence, of $St(y,K) \cap F$. It follows that $w$ is a vertex of $St(y,K)$. Since $w$ is not equal to $y$, this implies that $w$ and $y$ are joined by an edge of $K$. In other words, $w$ and $z$ are joined by an edge of $K$. 

In any case, $w$ and $z$ are joined by an edge of $K$. 

This shows that $\rho$ is a simplex of $K$. Since $\tau$ is a face of the simplex $\rho$ of $K$ and $y$ is a vertex of $\rho$, it follows that $\tau$ is a simplex of $St(y,K)$. 

This shows that $St(x,K) \subset St(y,K)$. Likewise, $St(y,K) \subset St(x,K)$. Hence, $St(x,K) = St(y,K)$. 

This proves that (\ref{xyanedge1}) implies (\ref{xyanedge2}). 

\end{proof}

\begin{prop} Let $K$ be a simplicial complex. Let $\mathcal{E}$ be a collection of exchangeable pairs of distinct vertices of $K$ with the property that no two distinct pairs in $\mathcal{E}$ have a common vertex. Then there exists a unique automorphism $\varphi_{\mathcal{E}} : K \rightarrow K$ such that (i) for each pair $\{x,y\}$ in $\mathcal{E}$, $\varphi_{\mathcal{E}}(x) = y$ and $\varphi_{\mathcal{E}}(y) = x$ and (ii) $\varphi_{\mathcal{E}}(z) = z$ for every vertex $z$ of $K$ which is not an element of some pair in $\mathcal{E}$. 
\label{prop:exchange} \end{prop}

\begin{proof} Let $\varphi : K_0 \rightarrow K_0$ be the unique involution which exchanges the two vertices in each pair in $\mathcal{E}$ and fixes every other vertex of $K$. Let $\tau$ be a simplex of $K$. We shall now show that $\varphi(\tau)$ is a simplex of $K$. Let $\tau_0$ be the set of all vertices $x$ of $\tau$ such that there does not exist a vertex $y$ of $K$ such that $\{x,y\} \in \mathcal{E}$;  
$\tau_1$ be the set of all vertices $x$ of $\tau$ such that there exists a vertex $y$ of $K$ such that $\{x,y\} \in \mathcal{E}$ and $\{x,y\} \cap \tau = \{x\}$; and 
$\tau_2$ be the set of all vertices $x$ of $\tau$ such that there exists a vertex $y$ of $K$ such that $\{x,y\} \in \mathcal{E}$ and $\{x,y\} \cap \tau = \{x,y\}$. 
Note that $\tau = \tau_0 \cup \tau_1 \cup \tau_2$. From the definition of $\varphi$, $\varphi(\tau_i) = \tau_i$, $i = 0,2$. 

Let $n$ be the number of elements of $\tau_1$. Suppose, on the one hand, that $n = 0$. Then, $\tau_1 = \emptyset$ and, hence, $\varphi(\tau) = \varphi(\tau_0 \cup \tau_2) = \varphi(\tau_0) \cup \varphi(\tau_2) = \tau_0 \cup \tau_2 = \tau$. Hence, $\varphi(\tau)$ is equal to the simplex $\tau$ of $K$. 

Suppose, on the other hand, that $n > 0$. Let $\tau_1 = \{x_j \vert 1 \leq j \leq n\}$. For each integer $j$ with $1 \leq j \leq n$, let $y_j$ be the unique vertex of $K$ such that $\{x_j,y_j\} \in \mathcal{E}$. From the definition of $\varphi$, $\varphi(\tau_1) = \{y_j \vert 1 \leq j \leq n\}$. It follows that $\varphi(\tau) = \tau_0 \cup \tau_2 \cup \{y_j \vert 1 \leq j \leq n\}$. 

Let $j$ be an integer with $1 \leq j \leq n$. Since $\{x_j,y_j\} \in \mathcal{E}$, $\{x_j,y_j\}$ is an exchangeable pair of vertices of $K$. Hence, there exists a simple exchange $\varphi_j : K \rightarrow K$ of $K$ exchanging $x_j$ and $y_j$. Since the distinct pairs $\{x_j,y_j\}$, $1 \leq j \leq n$, are in $\mathcal{E}$, they are disjoint. It follows that the composition $\varphi_1 \circ \ldots \circ \varphi_n : K \rightarrow K$ is an automorphism $\psi$ of $K$ such that $\psi(\tau_0) = \tau_0$, $\psi(\tau_2) = \tau_2$ and $\psi(x_j) = y_j$, $1 \leq j \leq n$. This implies that $\varphi(\tau) = \psi(\tau)$. Since $\psi : K \rightarrow K$ is an automorphism of $K$ and $\tau$ is a simplex of $K$, it follows that $\psi(\tau)$ is a simplex of $K$. That is to say, $\varphi(\tau)$ is a simplex of $K$. 

This shows that the involution $\varphi : K_0 \rightarrow K_0$ extends to a simplicial map $\varphi : K \rightarrow K$. Since $\varphi : K_0 \rightarrow K_0$ is an involution, its simplicial extension $\varphi : K \rightarrow K$ is also an involution and, hence, a simplicial automorphism of $K$. Hence, $\varphi : K \rightarrow K$ is a simplicial automorphism of $K$ such that (i) for each pair $\{x,y\}$ in $\mathcal{E}$, $\varphi_{\mathcal{E}}(x) = y$ and $\varphi_{\mathcal{E}}(y) = x$ and (ii) $\varphi_{\mathcal{E}}(z) = z$ for every vertex $z$ of $K$ which is not an element of some pair in $\mathcal{E}$. Since the stated conditions on $\varphi : K \rightarrow K$ determine the restriction $\varphi : K_0 \rightarrow K_0$, and any two simplicial maps which agree on the vertices of their common domain are equal, it follows that $\varphi : K \rightarrow K$ is the unique such automorphism of $K$.

\end{proof}

\begin{defn} Let $K$, $\mathcal{E}$, and $\varphi_{\mathcal{E}} : K \rightarrow K$ be as in Proposition \ref{prop:exchange}. We call the automorphism $\varphi_{\mathcal{E}} : K \rightarrow K$ of $K$ the {\it generalized exchange of $K$ associated to $\mathcal{E}$}. 
\label{defn:exchange} \end{defn}

If $\mathcal{F}$ and $\mathcal{G}$ are subsets of $\mathcal{E}$, then $\varphi_{\mathcal{F}} \circ \varphi_{\mathcal{G}} = \varphi_{\mathcal{F} \triangle \mathcal{G}}$. Hence, by Proposition \ref{prop:exchange}, we have the following result. 

\begin{prop} Let $K$ be a simplicial complex. Let $\mathcal{E}$ be a collection of exchangeable pairs of distinct vertices of $K$ with the property that no two distinct pairs in $\mathcal{E}$ have a common vertex. Then there exists a monomorphism $\Phi$ from the Boolean algebra $\mathcal{B}(\mathcal{E})$ of all subsets of $\mathcal{E}$ to $Aut(K)$ such that $\Phi(\mathcal{F}) = \varphi_{\mathcal{F}}$ for every subset $\mathcal{F}$ of $\mathcal{E}$. 
\label{prop:boolean} \end{prop}

\begin{defn} Let $K$ be a simplicial complex. Let $\mathcal{E}$ be a collection of exchangeable pairs of distinct vertices of $K$ with the property that no two distinct pairs in $\mathcal{E}$ have a common vertex.  The {\it Boolean subgroup of $Aut(K)$ corresponding to $\mathcal{E}$}, denoted by $B_{\mathcal{E}}$, is the image $\Phi(\mathcal{B}(\mathcal{E}))$ of the Boolean algebra $\mathcal{B}(\mathcal{E})$ under the monomorphism $\Phi$ of Proposition \ref{prop:boolean}. In particular, the Boolean subgroup $B_{\mathcal{E}}$ is naturally isomorphic to the Boolean algebra $\mathcal{B}(\mathcal{E})$. \label{defn:boolean} \end{defn}

\begin{prop} Let $K$ be a simplicial complex. Let $\mathcal{E}$ be a collection of exchangeable pairs of distinct vertices of $K$ with the property that no two distinct pairs in $\mathcal{E}$ have a common vertex. Let $\varphi \in Aut(K)$, $\mathcal{F} \subset \mathcal{E}$ and $\mathcal{G} = \varphi(\mathcal{F})$. Then $\mathcal{G}$ is a collection of exchangeable pairs of distinct vertices of $K$ with the property that no two distinct pairs in $\mathcal{G}$ have a common vertex.  Moreover, $\varphi \circ \Phi_{\mathcal{F}} \circ \varphi^{-1} = \Phi_{\mathcal{G}}$.
\label{prop:booleanconjugation} \end{prop}

\section{Complexes associated to $S$} 

In this section, we shall discuss some abstract simplicial complexes naturally associated to $S$. All of these complexes are finite-dimensional flag complexes.

The simplices of each of the complexes are finite collections of isotopy classes of subspaces of $S$ of a certain type which can be represented by disjoint subspaces of this type. The extended mapping class group $\Gamma^*(S)$ of $S$ acts naturally on each of these complexes via a natural action of the group of homeomorphisms of $S$ on the relevant subspaces.

For each of these complexes we say that an automorphism of the complex is {\it geometric} if it is induced by a homeomorphism of $S$. The main question addressed in this paper is whether, for each of these complexes, every automorphism is geometric. The affirmative answer to this question is known to hold for one of these complexes, the complex of curves, and a number of other complexes associated to $S$ (\cite{behrstockmargalit}, \cite{brendlemargalit}, (\cite{farbivanov}, \cite{irmak1}, \cite{irmak2}, \cite{irmak3}, \cite{irmakkorkmaz}, \cite{ivanov1}, \cite{korkmaz},\cite{luo}, \cite{margalit}, \cite{mv2}, \cite{shackleton}). See also \cite{mcc5} for a review.

\begin{defn} The {\it complex of curves of $S$}, $C(S)$, is the simplicial complex whose $n$-simplices are collections of $n + 1$ distinct isotopy classes of essential disjoint curves on $S$. \label{defn:curves} \end{defn}

Note that a finite collection of vertices of $C(S)$ forms a simplex of $C(S)$ if and only if each pair of vertices in this collection can be represented by disjoint curves on $S$. In other words, $C(S)$ is a flag complex. 

As shown in \cite{mm1}, \S2.2, $C(S)$ is empty when $S$ is a sphere with at most three holes; $C(S)$ is an infinite set of vertices when $S$ is a sphere with four holes or a torus with at most one hole; and $C(S)$ is connected when $S$ is not a sphere with at most four holes or a torus with at most one hole. 

\begin{defn} The {\it complex of domains of $S$}, $D(S)$, is the simplicial complex whose $n$-simplices are collections of $n + 1$ distinct isotopy classes of disjoint domains on $S$. \label{defn:complexofdomains} \end{defn}

Note that a finite collection of vertices of $D(S)$ forms a simplex of $D(S)$ if and only if each pair of distinct vertices in this collection can be represented by disjoint domains on $S$. In other words, $D(S)$ is a flag complex. 

Note that each system of domains on $S$ determines a simplex of $D(S)$ and each simplex of $D(S)$ is so determined.

Let $x, y \in D_0(S)$. On the one hand, if $x$ and $y$ are represented by disjoint domains $X$ and $Y$, then $y \in St(x,D(S))$. On the other hand, if $y \in St(x,D(S))$, it is not true, in general, that $x$ and $y$ are represented by disjoint domains on $S$. However, if $x \neq y$, then the following are equivalent: (i) $y \in St(x,D(S))$, (ii) $y \in Lk(x,D(S))$, (iii) $x$ and $y$ are represented by disjoint domains $X$ and $Y$. 

If $x$ is a vertex of $D(S)$ which is represented by an annulus, then we say that {\it $x$ is an annular vertex of $D(S)$}. 

\begin{prop} Let $x$ be an annular vertex of $D(S)$ and $y$ be a vertex of $D(S)$. Then the following are equivalent: 

\begin{enumerate} 

\item \label{starofannulus1} $y$ is a vertex of $St(x,D(S))$.

\item \label{starofannulus2} $\{x,y\}$ is a simplex of $D(S)$.

\item \label{starofannulus3} $x$ and $y$ are represented by disjoint domains $X$ and $Y$.

\end{enumerate} 

\label{prop:starofannulus} \end{prop}

\begin{defn} The {\it truncated complex of domains of $S$}, $D^2(S)$, is the induced subcomplex of $D(S)$ corresponding to those vertices of $D(S)$ which are not represented by biperipheral pairs of pants.\label{defn:truncated} \end{defn}

Note that a finite collection of vertices of $D^2(S)$ forms a simplex of $D^2(S)$ if and only if each pair of vertices in this collection can be represented by disjoint domains on $S$. In other words, $D^2(S)$ is a flag complex.

Note that $D^2(S) = D(S)$ when $b \leq 1$. In particular, $D^2(S) = D(S)$ for any closed surface $S$. 

A {\it biperipheral curve on $S$} is a curve on $S$ which is a boundary component of a biperipheral pair of pants.

There is a unique projection $\pi : D(S) \rightarrow D^2(S)$ which sends each vertex of $D^2(S)$ to itself and sends each remaining vertex of $D(S)$ to the vertex of $D^2(S)$ represented by a regular neighborhood of the unique essential boundary component of any biperipheral pair of pants representing this vertex.

Note that the for each vertex $x$ of $D^2(S)$ which is not represented by a regular neighborhood of a biperipheral curve on $S$ the fiber $\pi^{-1}(x)$ of $\pi: D(S) \rightarrow D^2(S)$ above $x$ is equal to $\{x\}$. 

Suppose that $x$ is a vertex of $D^2(S)$ which is represented by a regular neighborhood of a biperipheral curve $\gamma$ on $S$. 

Suppose that $S$ is a sphere with four holes. Then the fiber $\pi^{-1}(x)$ of $\pi: D(S) \rightarrow D^2(S)$ above $x$ is the triangle of $D(S)$ induced by the vertices of $D(S)$ corresponding to a regular neighborhood of $\gamma$ on $S$ and the two biperipheral pairs of pants on $S$ of which $\gamma$ is a boundary component. 

Suppose that $S$ is not a sphere with four holes. Then the fiber $\pi^{-1}(x)$ of $\pi: D(S) \rightarrow D^2(S)$ above $x$ is the edge of $D(S)$ induced by the vertices of $D(S)$ corresponding to a regular neighborhood of $\gamma$ on $S$ and the unique biperipheral pair of pants on $S$ of which $\gamma$ is a boundary component. 

Let $\alpha$ be an essential curve on $S$. Let $N_\alpha$ be a regular neighborhood of $\alpha$ on $S$. Note that $N_\alpha$ is an essential annulus on $S$. There is a natural inclusion $C(S) \hookrightarrow D^2(S)$ which maps the vertex of $C(S)$ represented by a curve $\alpha$ on $S$ to the vertex of $D^2(S)$ represented by the essential annulus $N_\alpha$ on $S$.  Note that the image $\iota(C(S))$ of $C(S)$ in $D^2(S)$ is the subcomplex of $D^2(S)$ induced by the set of annular vertices of $D^2(S)$. 

\section{Recognizing annular vertices in $D^2(S)$} 

\begin{prop} Suppose that $S$ is not a torus with one hole. Let $x \in D^2_0(S)$. Then the following are equivalent:

\begin{enumerate} 

\item \label{D2annuli1} $x$ is an annular vertex of $D^2(S)$.

\item \label{D2annuli2} For each vertex $y$ of $D^2(S)$ which is not equal to $x$, $St(x,D^2(S))$ is not contained in $St(y,D^2(S))$.

\end{enumerate}

\label{prop:D2annuli} \end{prop} 

\begin{proof} Since $D^2(S)$ is a flag complex, requiring property (\ref{D2annuli2}) of a vertex $x$ of $D^2(S)$ is equivalent to requiring that for each vertex $y$ of $D^2(S)$ which is not equal to $x$, there exists a vertex $z$ of $D^2(S)$ such that $\{x,z\}$ is a simplex of $D^2(S)$ and $\{y,z\}$ is not a simplex of $D^2(S)$. 

We begin by proving that (\ref{D2annuli1}) implies (\ref{D2annuli2}). To this end, let $x$ be an annular vertex of $D^2(S)$ and $y$ be a vertex of $D^2(S)$ such that $y \neq x$. 

We shall deduce (\ref{D2annuli2}) by contradiction. To this end, suppose that:

\begin{quote} (*) For every vertex $z$ of $D^2(S)$ such that $\{x,z\}$ is a simplex of $D^2(S)$, $\{y,z\}$ is a simplex of $D^2(S)$. 
\label{quote:star} \end{quote}

\noindent Since $x$ is an annular vertex of $D^2(S)$, there exists an essential curve $\alpha$ on $S$ such that $x$ is represented by regular neighborhoods of $\alpha$ on $S$.
 
Choose a maximal system $\mathcal{C}$ of curves on $S$ containing $\alpha$. 

Let $R$ be a regular neighborhood of the support $|\mathcal{C}|$ of $\mathcal{C}$ on $S$. For each curve $\beta$ in the system $\mathcal{C}$, let $R_\beta$ be the unique component of $R$ which contains $\beta$ and $x_\beta = [R_\beta] \in D^2(S)$.

Let $\beta \in \mathcal{C}$. Since $x_\alpha = x$, it follows that $\{x,x_\beta\}$ is a simplex of $D^2(S)$ and, hence, by condition (*), $\{y,x_\beta\}$ is a simplex of $D^2(S)$. 

In particular, $\{y,x\} = \{y,x_\alpha\}$ is a simplex of $D^2(S)$. Since $y$ is not equal to $x$, this implies that $\{y,x\}$ is an edge of $D^2(S)$. Hence, $S$ is neither a sphere with at most four holes nor a closed torus. 

In particular, the Euler characteristic of $S$ is negative. Hence, the maximal system $\mathcal{C}$ of curves on $S$ is a pants decomposition of $S$. Let $\mathcal{P}$ be the collection of pairs of pants of $\mathcal{C}$. 
 
Let $Y$ be a domain on $S$ representing $y$. 

Since $\{y,x_\beta\}$ is a simplex of $D^2(S)$ for every curve $\beta$ in the pants decomposition $\mathcal{C}$ of $S$, it follows that $Y$ is a domain on $S$ which is isotopic on $S$ either to $R_\beta$ for some $\beta$ in $\mathcal{C}$ or to $P$ for some pair of pants $P$ in $\mathcal{P}$. 

Suppose that $Y$ is isotopic on $S$ to $R_\beta$ for some $\beta$ in $\mathcal{C}$. 
Then $x_\beta = y \neq x = x_\alpha$ and, hence, $\beta \neq \alpha$. 

Since $\alpha$ and $\beta$ are disjoint nonisotopic essential curves on $S$, it follows from Proposition \ref{prop:detectisotopy} that there exists a curve $\gamma$ on $S$ such that $i(\alpha,\gamma) = 0$ and $i(\beta,\gamma) \neq 0$.

Let $Z$ be a regular neighborhood of $\gamma$ on $S$ and $z = [Z] \in D^2(S)$. Since $i(\alpha,\gamma) = 0$, it follows that $\{x,z\}$ is a simplex of $D^2(S)$. Hence, by condition (*), $\{y,z\}$ is a simplex of $D^2(S)$. Since $\{y,z\}$ is a simplex of $D(S)$ and $z$ is an annular vertex, it follows from Proposition \ref{prop:starofannulus} that $Y$ is isotopic on $S$ to a domain which is disjoint from $Z$. Hence, $i(\beta,\gamma) = 0$, which is a contradiction. 

Thus, $Y$ is isotopic on $S$ to some pair of pants $P$ in $\mathcal{P}$. 

Since $P$ represents the vertex $y$ of $D^2(S)$, $P$ is not a biperipheral pair of pants on $S$. Since $S$ is not a torus with one hole and $P$ is a domain on $S$ which is a nonbiperipheral pair of pants on $S$, it follows that there exists a pair of distinct nonisotopic essential boundary components, $\epsilon$ and $\eta$, of $P$. 

Since each essential boundary component of a pair of pants of a pants decomposition of $S$ is isotopic to one of the curves of the pants decomposition, there exist distinct curves, $\beta$ and $\delta$ of $\mathcal{C}$ such that $\epsilon$ and $\eta$ are isotopic on $S$ to $\beta$ and $\delta$. 

Since $\beta \neq \delta$, we may assume that $\alpha \neq \beta$. As before, it follows from Proposition \ref{prop:detectisotopy} that there exists a curve $\gamma$ on $S$ such that $i(\alpha,\gamma) = 0$ and $i(\beta,\gamma) \neq 0$.

Let $Z$ be a regular neighborhood of $\gamma$ on $S$ and $z = [Z] \in D^2(S)$. Since $i(\alpha,\gamma) = 0$, it follows that $\{x,z\}$ is a simplex of $D^2(S)$. Hence, by condition (*), $\{y,z\}$ is a simplex of $D^2(S)$. Since $\{y,z\}$ is a simplex of $D(S)$ and $z$ is an annular vertex, it follows from Proposition \ref{prop:starofannulus} that $Y$ is isotopic on $S$ to a domain which is disjoint from $Z$. Since $\beta$ is isotopic on $S$ to $\epsilon$ and $\epsilon \subset Y$, it follows that $i(\beta,\gamma) = i(\epsilon,\gamma) = 0$, which is a contradiction.

This shows that (\ref{D2annuli1}) implies (\ref{D2annuli2}).

We shall now show that (\ref{D2annuli2}) implies (\ref{D2annuli1}). To this end, suppose that the vertex $x$ of $D^2(S)$ is not an annular vertex of $D^2(S)$. 

Let $X$ be a domain on $S$ representing $x$. Since $x$ is not an annular vertex of $D^2(S)$, $X$ is not an annulus. 

Since $X$ is a domain on $S$, $X$ has an essential boundary component on $S$. Let $Y$ be a regular neighborhood of an essential boundary component of $X$ on $S$ and $y = [Y] \in D^2(S)$. Then $Y$ is isotopic on $S$ to a domain on $S$ which is disjoint from $X$. Since $X$ is not an annulus and $Y$ is an annulus, $Y$ is not isotopic to $X$ on $S$. It follows that $y$ is not equal to $x$ and $\{x,y\}$ is an edge of $D^2(S)$. 

Suppose that $z$ is a vertex of $D^2(S)$ such that $\{x,z\}$ is a simplex of $D^2(S)$. That is to say, suppose that either $x = z$ or $\{x,z\}$ is an edge of $D^2(S)$. 

Suppose, on the one hand, that $z = x$. Then $\{y,z\}$ is equal to the simplex $\{x,y\}$ of $D^2(S)$. 

Suppose, on the other hand, that $\{x,z\}$ is an edge of $D^2(S)$. Then $z$ is represented by a domain $Z$ on $S$ which is disjoint from the domain $X$ on $S$. Since $Y$ is a regular neighborhood of an essential boundary component of $X$, it follows that $Y$ is isotopic to a domain on $S$ which is disjoint from $Z$. 
This implies that if $Z$ is isotopic to $Y$, then $\{y,z\}$ is equal to the simplex $\{y\}$ of $D^2(S)$, whereas, if $Z$ is not isotopic to $Y$ on $S$, then $\{y,z\}$ is an edge of $D^2(S)$. In any case, $\{y,z\}$ is a simplex of $D^2(S)$. 

This shows that (\ref{D2annuli2}) implies (\ref{D2annuli1}).

\end{proof}

\begin{cor} Suppose that $S$ is not a torus with one hole. Then every simplicial automorphism of $D^2(S)$ restricts to a simplicial automorphism of the subcomplex $\iota(C(S))$ of $D^2(S)$ induced from the set of annular vertices of $D^2(S)$.
\label{cor:D2annuli} \end{cor}

\begin{proof} Let $\varphi \in Aut(D^2(S))$. 

Suppose that $x$ is an annular vertex of $D^2(S)$ (i.e. a vertex of $\iota(C(S))$). By Proposition \ref{prop:D2annuli}, for each vertex $y$ of $D^2(S)$ which is not equal to $x$, there exists a vertex $z$ of $D^2(S)$ such that $\{x,z\}$ is a simplex of $D^2(S)$ and $\{y,z\}$ is not a simplex of $D^2(S)$.

Let $u = \varphi(x)$. Since $\varphi \in Aut(D^2(S))$, it follows that for each vertex $v$ of $D^2(S)$ which is not equal to $u$, there exists a vertex $w$ of $D^2(S)$ such that $\{u,w\}$ is a simplex of $D^2(S)$ and $\{v,w\}$ is not a simplex of $D^2(S)$. Hence, by Proposition \ref{prop:D2annuli}, $u$ is a vertex of $\iota(C(S))$. This shows that $\varphi$ maps the zero skeleton of $\iota(C(S))$ into the zero skeleton of $\iota(C(S))$. 

Note that a simplex $\sigma$ of $D^2(S)$ is a simplex of $\iota(C(S))$ if and only if each of its vertices is a vertex of $\iota(C(S))$. It follows that $\varphi$ restricts to a simplicial map $\mu : \iota(C(S)) \rightarrow \iota(C(S))$. Likewise, the simplicial automorphism $\varphi^{-1} : D^2(S) \rightarrow D^2(S)$ restricts to a simplicial map $\lambda : \iota(C(S)) \rightarrow \iota(C(S))$. 

Note that the restrictions $\mu : \iota(C(S)) \rightarrow \iota(C(S))$ and $\lambda : \iota(C(S)) \rightarrow \iota(C(S))$ of $\varphi : D^2(S) \rightarrow D^2(S)$ and $\varphi^{-1} : D^2(S) \rightarrow D^2(S)$ are inverse simpicial maps. Hence, $\mu : \iota(C(S)) \rightarrow \iota(C(S))$ is a simplicial isomorphism. This proves that $\varphi : D^2(S) \rightarrow D^2(S)$ restricts to a simplicial automorphism $\varphi : \iota(C(S)) \rightarrow \iota(C(S))$.

\end{proof}

\section{Distinguishing vertices of $D(S)$ via their annular links}

\begin{defn} Let $x$ be a vertex of $D(S)$. The {\it annular link of $x$ in $D(S)$} is the subcomplex $Ann(x)$ of $D(S)$ consisting of those simplices of $Lk(x,D(S))$ all of whose vertices are annular.  
\label{defn:annx} \end{defn}

\begin{prop} Suppose that $S$ is neither a sphere with four holes nor a torus with at most one hole. Let $x$ and $y$ be vertices of $D(S)$. Suppose that $x$ is annular. Then the following are equivalent: 

\begin{enumerate} 

\item \label{nestedann3.1} $Ann(x) \subset Ann(y)$

\item \label{nestedann3.2} $x = y$ or there exist disjoint domains $X$ and $Y$ on $S$ representing $x$ and $y$ such that $X$ is an annulus on $S$, $Y$ is a biperipheral pair of pants on $S$, and $X \cup Y$ has exactly two codomains, exactly one of which is an annulus joining $X$ to $Y$. 

\end{enumerate} 

\label{prop:nestedann3} \end{prop} 

\begin{proof}  We begin by proving that (\ref{nestedann3.1}) implies (\ref{nestedann3.2}). Suppose that $Ann(x) \subset Ann(y)$. 

Since $x$ is annular, $x$ is represented by a regular neighborhood $X$ of an essential curve $\alpha$ on $X$. 

Since $S$ has an essential curve $\alpha$, $S$ is not a sphere with at most three holes. Since $S$ is also not a closed torus, there exists a pants decomposition $\mathcal{C}$ of $S$ containing $\alpha$. Let $R$ be a regular neighborhood of the support $|\mathcal{C}|$ of $\mathcal{C}$ on $S$ and $\mathcal{P}$ be the collection of codomains of $R$ on $S$. We may assume that $X$ is the unique component of $R$ which contains the element $\alpha$ of $\mathcal{C}$.

Note that each element of $\mathcal{P}$ is a pair of pants on $S$. 

Let $\beta$ be an element of $\mathcal{C}$ which is not equal to $\alpha$. Then a regular neighborhood $Z$ of $\beta$ on $S$ represents a vertex $z$ of $Ann(x)$ and, hence, of $Ann(y)$. It follows that $y$ is represented by a domain $Y$ on $S$ which is disjoint from and not isotopic to each of the components of a regular neighborhood $W$ of the union of all the elements of $\mathcal{C}$ which are not equal to $\alpha$. Since $Y$ is connected, $Y$ is contained in a codomain of $W$ on $S$.

Note that the unique codomain $V$ of $W$ on $S$ which contains $X$ is equal to the union of $X$ with those elements of $\mathcal{P}$ which share at least one common essential boundary component with $X$. If there is exactly one such element of $\mathcal{P}$, then $W$ is a torus with one hole. Otherwise, there are exactly two such elements of $\mathcal{P}$ and $W$ is a sphere with four holes. 

Every other codomain $U$ of $W$ on $S$ is a pair of pants on $S$ all of whose essential boundary components are isotopic to elements of $\mathcal{C}$ which are not equal to $\alpha$. 

Suppose that $Y$ is contained in one of these other codomains $U$ of $W$ on $S$. Since $U$ is a pair of pants and $Y$ is a domain on $S$ contained in $U$, $Y$ is isotopic on $S$ to a domain $Y_1$ on $S$ such that $Y_1$ is equal to $U$ or $Y_1$ is a regular neighborhood of an essential boundary component of $U$ on $S$. Note, in any case, that an essential boundary component $\delta$ of $U$ is contained in $Y_1$.   

By assumption, $\delta$ is isotopic to an element $\beta$ of $\mathcal{C}$ which is not equal to $\alpha$. Since $\alpha$ and $\beta$ are distinct elements of the pants decomposition $\mathcal{C}$, $\alpha$ and $\beta$ are disjoint nonisotopic essential curves on $S$. It follows from Proposition \ref{prop:detectisotopy} that there exists a curve $\gamma$ on $S$ such that $i(\gamma,\alpha) = 0$ and $i(\gamma,\beta) \neq 0$. Since $\delta$ is isotopic to $\beta$, $i(\gamma,\delta) = i(\gamma,\beta)$. Hence, $i(\gamma,\delta) \neq 0$.

It follows that a regular neighborhood $Z$ of $\gamma$ on $S$ represents a vertex $z$ of $Ann(x)$ and, hence, of $Ann(y)$. Since $Y_1$ represents the vertex $y$ of $D(S)$ and $Z$ represents the vertex $z$ of $D(S)$, it follows that $Z$ is isotopic on $S$ to a domain $Z_1$ on $S$ which is disjoint from $Y_1$. Thus, $\gamma$ is isotopic on $S$ to a curve $\gamma_1$ on $S$ which is disjoint from $\delta$. Thus, $i(\gamma,\delta) = i(\gamma_1,\delta) = 0$, which is a contradiction. 

Hence, $Y$ is not contained in one of these other codomains $U$ of $W$ on $S$. It follows that $Y$ is not isotopic on $S$ to a domain which is contained in one of the other codomains $U$ of $W$ on $S$.

It follows that $Y$ is contained in the unique codomain $V$ of $W$ on $S$ which contains $X$. 

Since $S$ is not a sphere with four holes nor a torus with one hole, $V$ has an essential boundary component $\delta$ on $S$. Note that $\delta$ is isotopic on $S$ to an element $\beta$ of $\mathcal{C}$ which is not equal to $\alpha$. 

As before, it follows from Proposition \ref{prop:detectisotopy} that there exists a curve $\gamma$ on $S$ such that $i(\gamma,\alpha) = 0$ and $i(\gamma,\beta) \neq 0$. Since $i(\gamma,\alpha) = 0$, we may assume that $\gamma$ is disjoint from $\alpha$. Since $\delta$ is isotopic to $\beta$, $i(\gamma,\delta) = i(\gamma,\beta)$. Hence, $i(\gamma,\delta) \neq 0$.

It follows that a regular neighborhood $Z$ of $\gamma$ on $S$ represents a vertex $z$ of $Ann(x)$ and, hence, of $Ann(y)$. Hence, $Y$ is isotopic on $X$ to a domain $Y_1$ on $X$ which is disjoint from $Z$ and, hence, from $\gamma$. Note that $X$ and $Y_1$ are both domains on $S$ which are contained in $V$ and are disjoint from $\gamma$. 

We may assume that the number of points of intersection of $\gamma$ with each essential boundary component $\epsilon$ of $V$ is equal $i(\gamma,\epsilon)$. Since $i(\gamma,\delta) \neq 0$, it follows that $\gamma \cap V$ is a nonempty disjoint union of properly embedded essential arcs on $V$. 

Suppose, on the one hand, that $V$ is a torus with one hole. Then, since $X$ and $Y_1$ are both domains on $S$ contained in $V$ and disjoint from a properly embedded essential arc on $V$, it follows that $X$ and $Y_1$ are isotopic annuli on $S$ and, hence, $x = y$. 

Suppose, on the other hand, that $V$ is a sphere with four holes. Then, since $X$ and $Y_1$ are both domains on $S$ contained in $V$ and disjoint from a properly embedded essential arc on $V$, it follows that $Y_1$ is isotopic to a domain $Y_2$ on $S$ which is contained in one of the two elements $P$ of $\mathcal{P}$ which share an essential boundary component with the annulus $X$. 

Suppose that $Y_2$ is isotopic to a regular neighborhood of an essential boundary component $\sigma$ of $P$. Note that $\sigma$ is isotopic to an element $\beta$ of $\mathcal{P}$. Since $Y_2$ is not isotopic to a regular neighborhood of any element of $\mathcal{C}$ which is not equal to $\alpha$, it follows that $\beta$ is equal to $\alpha$. This implies that $X$ is isotopic to $Y_2$ on $S$ and, hence, $x = y$. 

Hence, we may assume that $Y_2$ is not isotopic to a regular neighborhood of any essential boundary component of $P$. Since $Y_2$ is a domain on $S$ contained in the pair of pants $P$, it follows that $Y_2$ is isotopic to $P$ on $S$. 

Suppose that there exists an essential boundary component $\tau$ of $P$ such that $\tau$ is not isotopic to $\alpha$ on $S$. Then $\tau$ is isotopic to an element $\beta$ of $\mathcal{C}$ which is not equal to $\alpha$. 

As before, it follows from Proposition \ref{prop:detectisotopy} that there exists a curve $\gamma$ on $S$ such that $i(\gamma,\alpha) = 0$ and $i(\gamma,\beta) \neq 0$. 

It follows that a regular neighborhood $Z$ of $\gamma$ on $S$ represents a vertex $z$ of $Ann(x)$ and, hence, of $Ann(y)$. 

Since $P$ represents $y$, it follows that $Z$ is isotopic to a domain on $S$ which is disjoint from $P$. Hence, $\gamma$ is isotopic to a curve $\gamma_1$ on $S$ which is disjoint from $\tau$. 

This implies that $i(\gamma,\beta) = i(\gamma_1,\tau) = 0$, which is a contradiction. 

Hence, each essential boundary component of $P$ is isotopic to $\alpha$ on $S$.

It follows that $P$ is either a monoperiperipheral pair of pants sharing both of its essential boundary components with $X$ or a biperipheral pair of pants sharing its unique essential boundary component with $X$. In the former case, it follows that $S$ is a torus with one hole, which is a contradiction. Hence, the latter case holds. 

Since $Y_2$ is a nonannular domain on $S$ contained in the biperipheral pair of pants $P$ on $S$, it follows that $Y_2$ is a biperipheral pair of pants on $S$ whose unique codomain on $P$ is an annulus on $S$ joining $X$ to $Y_2$. 

This completes the proof that (\ref{nestedann3.1}) implies (\ref{nestedann3.2}). It remains to prove that (\ref{nestedann3.2}) implies (\ref{nestedann3.1}).

If $x = y$, then $Ann(x) = Ann(y)$ and, hence, $Ann(x) \subset Ann(y)$. 

Suppose that $x$ and $y$ are represented by disjoint domains $X$ and $Y$ on $S$ such that $X$ is an annulus on $S$, $Y$ is a biperipheral pair of pants on $S$, and $X \cup Y$ has exactly two codomains, exactly one of which is an annulus joining $X$ to $Y$. 

Let $P$ be the unique codomain of $X$ on $S$ such that $Y$ is contained in $P$. Note that $P$ is a biperipheral pair of pants on $S$.

Suppose that $z$ is an element of $Ann(x)$. Then $z$ is represented by an annulus $Z$ on $S$ which is disjoint from and not isotopic to $X$. 

Since $Z$ is connected and disjoint from $X$, $Z$ is contained in a codomain $Q$ of $X$ on $S$. 

Suppose that $Q$ is equal to $P$. Then since $Z$ is an annular domain on $S$ and $P$ is a biperipheral pair of pants on $S$, it follows that $Z$ is isotopic to a regular neighborhood of the unique essential boundary component of $P$ on $S$. Since every essential boundary component of a codomain of $X$ on $S$ is an essential boundary component of the annulus $X$, it follows that $Z$ is isotopic to $X$ on $S$, which is a contradiction. Hence, $Q$ is not equal to $P$. 

Since any two distinct codomains of $X$ on $S$ are disjoint, it follows that $Z$ is disjoint from $P$ and, hence, from $Y$. Note that the annulus $Z$ is not isotopic on $S$ to the pair of pants $Y$. Hence, the vertex $z$ of $D(S)$ represented by the annulus $Z$ is an element of $Ann(y)$. 

This proves that (\ref{nestedann3.2}) implies (\ref{nestedann3.1}), completing the proof. 

\end{proof}

\begin{prop} Suppose that $S$ is neither a sphere with four holes nor a torus with at most one hole. Let $x$ and $y$ be vertices of $D(S)$. Suppose that $x$ is annular. Then $Ann(x) = Ann(y)$ if and only if $x = y$.
\label{prop:ann2} \end{prop} 

\begin{proof} It suffices to prove that $Ann(x) = Ann(y)$ implies $x = y$. 

To this end, suppose that $Ann(x) = Ann(y)$. Then $x$ is annular and $Ann(x) \subset Ann(y)$. It follows from Proposition \ref{prop:nestedann3} that either (i) $x = y$ or (ii) $x$ and $y$ are represented by an annulus $X$ on $S$ and a biperipheral pair of pants $Y$ on $S$ such that $X \cup Y$ has exactly two codomains, exactly one of which is an annulus joining $X$ to $Y$. 

Suppose that (ii) holds. Then $x$ is a vertex of $Ann(y)$. That is to say, since $Ann(x) = Ann(y)$, $x$ is a vertex of $Ann(x)$. Since $Ann(x)$ is a subcomplex of $Lk(x,D(S))$, it follows that $x$ is a vertex of $Lk(x,D(S))$, which is a contradiction. Hence, (ii) does not hold.

It follows that (i) holds. That is to say, it follows that $x = y$, completing the proof. 
 
\end{proof}

\begin{prop} Let $x$ and $y$ be vertices of $D(S)$. Suppose that $\{x,y\}$ is a simplex of $D(S)$. Then $Ann(x) \subset Ann(y)$ if and only if either $x = y$ or $x$ and $y$ are represented by disjoint domains $X$ and $Y$ on $S$ such that $Y$ is a pair of pants with each of its essential boundary components on $S$ joined to $X$ by annuli. 
\label{prop:nestedann4} \end{prop} 

\begin{proof} Suppose, on the one hand, that $Ann(x) \subset Ann(y)$. 

We may assume that $x$ is not equal to $y$. Then, since $\{x,y\}$ is a simplex of $D(S)$, $\{x,y\}$ is an edge of $D(S)$. It follows that $x$ and $y$ are represented by disjoint domains $X$ and $Y$ on $S$ which are not isotopic to one another on $S$. 

Suppose that $Y$ is a nonelementary domain on $S$. It follows from Proposition \ref{prop:nonelemdom} that there exist curves $\alpha$ and $\beta$ on $S$ such that $i(\alpha,\beta) \neq 0$ and $\alpha$ and $\beta$ are contained in the interior of $Y$. 

Let $U$ and $V$ be regular neighborhoods of $\alpha$ and $\beta$ in the interior of $Y$. Suppose that $V$ is isotopic to $X$ on $S$. Then $\beta$ is isotopic on $S$ to a curve $\beta_1$ which is contained in the interior of $X$ and, hence, is disjoint from $Y$. Since $\beta$ is isotopic to $\beta_1$ on $S$, $i(\alpha,\beta) = i(\alpha,\beta_1)$. Since $\alpha$ is contained in $Y$ and $\beta_1$ is disjoint from $Y$, it follows that $\alpha$ and $\beta_1$ are disjoint and, hence, $i(\alpha,\beta_1) = 0$. We conclude that $i(\alpha,\beta) = 0$ which is a contradiction. 

ence, $V$ is not isotopic to $X$ on $S$. Since $V$ is contained in $Y$ and $X$ and $Y$ are disjoint, $X$ and $V$ are disjoint domains on $S$. It follows that $V$ represents an annular vertex $v$ of $Lk(x,D(S))$. This implies that $v$ is a vertex of $Ann(x)$ and, hence, of $Ann(y)$. 

Since $V$ represents $v$ and $Y$ represents $y$, it follows that $V$ is isotopic on $S$ to a domain on $S$ which is disjoint from $Y$. Since $\beta$ is contained in $V$, it follows that $\beta$ is isotopic on $S$ to a curve $\beta_2$ which is disjoint from $Y$. Again, this implies that $i(\alpha,\beta) = i(\alpha,\beta_2) = 0$, which is a contradiction. 

It follows that $Y$ is an elementary domain on $S$. 

Suppose that $Y$ is an annulus. Since $X$ and $Y$ are disjoint nonisotopic domains on $S$, it follows that the annular vertex $y$ of $D(S)$ represented by $Y$ is a vertex of $Ann(x)$ and, hence of $Ann(y)$. Since $Ann(y)$ is a subcomplex of $Lk(y,D(S))$, it follows that $y$ is a vertex of $Lk(y,D(S))$, which is a contradiction. 

Hence, $Y$ is not an annulus. Since $Y$ is an elementary domain on $S$, it follows that $Y$ is a pair of pants. 

Let $\beta$ be an essential boundary component of $Y$ on $S$. Suppose that $\beta$ is not isotopic to any essential boundary component of $X$ on $S$. Then, by Proposition \ref{prop:detectisotopy}, there exists an essential curve $\gamma$ on $S$ such that $i(\gamma,\alpha) = 0$ for every essential boundary component $\alpha$ of $X$ on $S$ and $i(\gamma,\beta) \neq 0$. 

Since $Y$ is disjoint from $X$, it follows that a regular neighborhood $W$ of $\gamma$ on $S$ represents a vertex $w$ of $Ann(x)$ and, hence, of $Ann(y)$. It follows that $\gamma$ is isotopic on $S$ to a curve $\gamma_1$ on $S$ which is disjoint from $Y$. It follows that $i(\gamma,\beta) = i(\gamma_1,\beta) = 0$, which is a contradiction. 

Hence, the essential boundary component $\beta$ of $Y$ on $S$ is isotopic on $S$ to some essential boundary component $\alpha$ of $X$ on $S$. Since $X$ and $Y$ are disjoint, it follows that there is an annulus $A$ on $S$ whose boundary components are $\alpha$ and $\beta$. 

Since $Y$ is a pair of pants, it follows that $A \cap Y = \beta$. 
Moreover, either $A \cap X = \alpha$ or $X \subset A$. In the former case, $A$ is an annulus joining $\beta$ to $X$. 

Suppose $X \subset A$. Then $X$ is an annulus contained in $A$. It follows that $A = X \cup B$, where $B$ is an annulus joining $\beta$ to $X$. 

In any case, $\beta$ is joined to $X$ by an annulus. 

This proves the ``only if'' direction. It remains to prove the ``if'' direction.

If $x = y$, then $Ann(x) = Ann(y)$ and, hence, $Ann(x) \subset Ann(y)$. 

Suppose that $x$ and $y$ are represented by disjoint domains $X$ and $Y$ on $S$ such that $Y$ is a pair of pants with each of its essential boundary components on $S$ joined to $X$ by annuli.

Since $Y$ is disjoint from $X$ on $S$ and each of the essential boundary components of $Y$ on $S$ is joined to $X$ by an annulus, it follows that $Y$ is isotopic on $S$ to the unique codomain $Y_1$ of $X$ on $S$ which contains $Y$.

Let $z$ be an element of $Ann(x)$. Then $z$ is represented by an annulus $Z$ on $S$ which is disjoint from $X$ and, hence, is contained in a codomain $W$ of $X$ on $S$. 

Suppose that $W$ is not equal to $Y_1$. Then $Z$ is disjoint from $Y_1$ and, hence, from $Y$. Since $Y$ is a pair of pants, it follows that the annulus $Z$ is not isotopic to $Y$ on $S$. Hence, the annular vertex $z$ of $D(S)$ represented by $Z$ is a vertex of $Ann(y)$. 

Suppose that $W$ is equal to $Y_1$. 

Since $Z$ is an annular domain on $S$ contained in the pair of pants $Y_1$ on $S$, it follows that $Z$ is isotopic on $S$ to an annulus $Z_1$ on $S$ which is disjoint from $Y_1$ and, hence, from $Y$. Note that the annulus $Z_1$ is not isotopic to the pair of pants $Y$ on $S$. Hence, the vertex $z$ of $D(S)$ represented by $Z_1$ is a vertex of $Ann(y)$. 

In any case, $z$ is an element of $Ann(y)$.

Again, we conclude that $Ann(x) \subset Ann(y)$. 

In any case, $Ann(x) \subset Ann(y)$. 

This proves the ``if'' direction, completing the proof. 

\end{proof}

\begin{prop} Suppose that $S$ is neither a sphere with four holes nor a torus with at most one hole. Let $x$ and $y$ be vertices of $D(S)$. Suppose that $\{x,y\}$ is not a simplex of $D(S)$. Then $Ann(x) \subset Ann(y)$ if and only if $x$ and $y$ are represented by domains $X$ and $Y$ on $S$ such that $Y$ is a domain on $X$. 
\label{prop:nestedann5} \end{prop} 

\begin{proof} Since $\{x,y\}$ is not a simplex of $D(S)$, $x \neq y$. 

Suppose that $Ann(x) \subset Ann(y)$. 

Suppose that $x$ is an annular vertex of $D(S)$. Since $S$ is neither a sphere with four holes nor a torus with at most one hole and $x \neq y$, it follows from Proposition \ref{prop:nestedann3} that $x$ and $y$ are represented by disjoint domains $X$ and $Y$ on $S$. Hence, $\{x,y\}$ is a simplex of $D(S)$, which is a contradiction. 
Therefore, $x$ is not an annular vertex of $D(S)$. 

It follows that a regular neighborhood $Z$ of any essential boundary component $\alpha$ of $X$ represents a vertex $z$ of $Ann(x)$ and, hence, of $Ann(y)$. It follows that $y$ is represented by a domain $Y$ on $S$ which is not isotopic to a regular neighborhood of any essential boundary component of $X$ on $S$ and is disjoint from a regular neighborhood of the union of the essential boundary components of $X$ on $S$. 

Since $Y$ is disjoint from a regular neighborhood of the union of the essential boundary components of $X$ on $S$, either $Y$ is disjoint from $X$ or $Y$ is contained in $X$. Since $\{x,y\}$ is not a simplex of $D(S)$, $Y$ is not disjoint from $X$. Hence, $Y$ is contained in $X$. 

Since $Y$ is a domain on $S$ contained in the domain $X$ on $S$, it follows from Proposition \ref{prop:weakconv}, that $Y$ is isotopic on $S$ to either $X$, or a domain on $X$, or a regular neighborhood of an essential boundary component of $X$.

Since $x \neq y$, $Y$ is not isotopic on $S$ to $X$. Since $Y$ is not isotopic on $S$ to a regular neighborhood of any essential boundary component of $X$ on $S$, we conclude that $Y$ is isotopic to a domain $Y_1$ on $X$. 

Hence, $x$ and $y$ are represented by domains $X$ and $Y_1$ on $S$ such that $Y_1$ is a domain on $X$. 

This proves the ``only if'' direction. It remains to prove the ``if'' direction. 

Suppose that $x$ and $y$ are represented by domains $X$ and $Y$ on $S$ such that $Y$ is a domain on $X$.

Let $z$ be a vertex of $Ann(x)$. Then $z$ is represented by an annulus $Z$ on $S$ which is disjoint from $X$. Since $Y$ is contained in $X$, it follows that $Z$ is disjoint from $Y$. 

Suppose that $Y$ is isotopic to $Z$ on $S$. Then $Y$ is an annulus on $S$. Hence, $Y$ is a regular neighborhood of an essential curve $\alpha$ on $S$. Since $Y$ is a domain on $X$, $\alpha$ is an essential curve on $X$. 

It follows from Proposition \ref{prop:geomint}, that there exists an essential curve $\beta$ on $X$ such that the geometric intersection number of $\alpha$ and $\beta$ on $S$ is not equal to zero.

Since $Z$ is an annular domain on $S$, $Z$ is a regular neighborhood of an essential curve $\gamma$ on $S$. Since $\beta$ is contained in $X$ and $Z$ is disjoint from $X$, it follows that $\beta$ is disjoint from $\gamma$ and, hence,  $i(\gamma,\beta) = 0$. Since $i(\alpha,\beta) \neq 0$, it follows that $\alpha$ is not isotopic to $\gamma$ on $S$. This implies that $Y$ is not isotopic to $Z$. Since $Z$ is an annulus disjoint from $Y$ and not isotopic to $Y$ on $S$, it follows that the vertex $z$ of $D(S)$ represented by $Z$ is a vertex of $Ann(y)$. 

This proves the ``if'' direction, completing the proof. 

\end{proof}

\begin{prop} Suppose that $S$ is neither a sphere with four holes nor a torus with at most one hole. Let $x$ and $y$ be vertices of $D(S)$. Then the following are equivalent: 

\begin{enumerate} 

\item \label{ann3.1} $Ann(x) = Ann(y)$

\item \label{ann3.2} $x = y$ or there exist disjoint domains $X$ and $Y$ on $S$,
representing $x$ and $y$, which belong to one of the following cases:  

\begin{enumerate} 

\item \label{ann3a} $S$ is a torus with two holes, $X$ and $Y$ are monoperipheral pairs of pants on $S$, and $X \cup Y$ has exactly two codomains, both of which are annuli joining $X$ to $Y$.

\item \label{ann3b} $S$ is a closed surface of genus two, $X$ and $Y$ are pairs of pants on $S$, and $X \cup Y$ has exactly three codomains, all of which are annuli joining $X$ to $Y$.

\end{enumerate} 

\end{enumerate} 

\label{prop:ann3} \end{prop} 

\begin{proof} 

We begin by proving the ``only if'' direction. Suppose that $Ann(x) = Ann(y)$. 

If $x$ is annular, then, since $Ann(x) = Ann(y)$, it follows from Proposition \ref{prop:ann2} that $x = y$. Likewise, if $y$ is annular, then, since $Ann(y) = Ann(x)$, it follows from Proposition \ref{prop:ann2} that $y = x$. Hence, if either $x$ or $y$ is annular, then $x = y$. 

Thus, we may assume that neither $x$ nor $y$ is annular. 

We may assume that $x \neq y$. 

Suppose that $\{x,y\}$ is not a simplex. Then, since $Ann(x) \subset Ann(y)$, it follows from Proposition \ref{prop:nestedann5}, that $x$ and $y$ are represented by domains $X$ and $Y$ on $S$ such that $Y$ is a domain on $X$. Likewise, since $Ann(y) \subset Ann(x)$, it follows from Proposition \ref{prop:nestedann5}, that $y$ and $x$ are represented by domains $Y_1$ and $X_1$ on $S$ such that $Y_1$ is a domain on $X_1$. Thus $X$ is isotopic to a domain on $Y$ and $Y$ is isotopic to a domain on $X$, which contradicts Proposition \ref{prop:partialorder}.

Hence, $\{x,y\}$ is a simplex of $D(S)$. Since $Ann(x) \subset Ann(y)$ and $x \neq y$, it follows from Proposition \ref{prop:nestedann4} that $x$ and $y$ are represented by disjoint domains $X$ and $Y$ on $S$ such that $Y$ is a pair of pants with each of its essential boundary components on $S$ joined to $X$ by annuli. This implies that the number of essential boundary components of $Y$ on $S$ is less than or equal to the number of essential boundary components of $X$ on $S$. 

Likewise, since $Ann(y) \subset Ann(x)$ and $y \neq x$, it follows from Proposition \ref{prop:nestedann4} that $y$ and $x$ are represented by disjoint domains $Y_1$ and $X_1$ on $S$ such that $X_1$ is a pair of pants with each of its essential boundary components on $S$ joined to $Y_1$ by annuli. Again, this implies that the number of essential boundary components of $X_1$ on $S$ is less than or equal to the number of essential boundary components of $Y_1$ on $S$. 

Since $X$ and $X_1$ both represent $x$, $X$ is isotopic to the pair of pants $X_1$ on $S$. This implies that $X$ is a pair of pants on $S$ with the same number of essential boundary components on $S$ as $X_1$. Likwise, $Y_1$ is a pair of pants on $S$ with the same number of essential boundary components on $S$ as $Y$. Since the number of essential boundary components of $X_1$ on $S$ is less than or equal to the number of essential boundary components of $Y_1$ on $S$, it follows that the number of essential boundary components of $X$ on $S$ is less than or equal to the number of essential boundary components of $Y$ on $S$. Since the number of essential boundary components of $Y$ on $S$ is less than or equal to the number of essential boundary components of $X$ on $S$, we conclude that $X$ is a pair of pants on $S$ with the same number of essential boundary components on $S$ as the pair of pants $Y$ on $S$. 

Thus, $X$ and $Y$ are disjoint pairs of pants on $S$ with the same number $n$ of essential boundary components on $S$ and the $n$ essential boundary components of $X$ on $S$ are joined by disjoint annuli to the $n$ essential boundary components of $Y$ on $S$. 

Since $X$ is a pair of pants domain on $S$, $1 \leq n \leq 3$. If $n = 1$, then $S$ is a sphere with four holes, which is a contradiction. Hence, $2 \leq n \leq 3$. If $n = 2$, then $X$ and $Y$ satisfy case (\ref{ann3a}). If $n = 3$, then $X$ and $Y$ satisfy case (\ref{ann3b}).

This completes the proof of the ``only if'' direction. It remains to prove the ``if'' direction.

If $x = y$, then $Ann(x) = Ann(y)$. 

Suppose that $X$ and $Y$ are as in case (\ref{ann3a}). Note that the two codomains of $X \cup Y$ on $S$ are annuli which are disjoint from and not isotopic on $S$ to the pairs of pants $X$ and $Y$ on $S$. Hence, they represent vertices of $Ann(x)$ and $Ann(y)$. 

Suppose that $z$ is a vertex of $Ann(x)$. Then $z$ is represented by an annulus on $S$ which is contained in the unique codomain $Y_1$ of $X$ on $S$. Note that $Y_1$ is a pair of pants on $S$ which is isotopic to $Y$ on $S$. It follows that $Z$ is isotopic to a regular neighorhood of an essential boundary component $\alpha$ of $Y_1$ on $S$. Since $Y_1$ is a codomain of $X$ on $S$, $\alpha$ is an essential boundary component of $X$ on $S$. It follows that $Z$ is isotopic to one of the two codomains of $X \cup Y$ on $S$. This proves that $Ann(x)$ is the edge of $D(S)$ whose vertices are represented by the two codomains of $X \cup Y$ on $S$. 
Likewise, $Ann(y)$ is equal to this edge and, hence, $Ann(x) = Ann(y)$. 

Similarly, if $X$ and $Y$ are as in case (\ref{ann3b}), then $Ann(x) = Ann(y)$. 

In any case, $Ann(x) = Ann(y)$. 

This proves the ``if'' direction, completing the proof. 

\end{proof} 

\section{Automorphisms of $D^2(S)$ are geometric}

In this section, we prove that if $S$ is not a sphere with at most four holes, a torus with at most two holes, or a closed surface of genus two, then each automorphism of $D^2(S)$ is induced by a self-homeo\-morphism of $S$ which is uniquely defined up to isotopy on $S$. This will imply that, under the same hypothesis on $S$, $Aut(D^2(S)) \simeq \Gamma^*(S)$ and, if $b \leq 1$, $Aut(D(S)) \simeq \Gamma^*(S)$.
 
%%%%\begin{lem} Suppose that $S$ is neither a sphere with at most three holes nor a %%%%torus with at most one hole. Let $\varphi : D^2(S) \rightarrow D^2(S)$ be an %%%%automorphism of $D^2(S)$. Then there exists an automorphism $\tau : C(S) %%%%\rightarrow C(S)$ such that $\varphi \circ i = i \circ \tau$. \label{lem:cinv} %%%%\end{lem}

\begin{lem} Suppose that $S$ is not a torus with one hole. Let $i : C(S) \rightarrow D^2(S)$ be the natural inclusion corresponding to forming regular neighborhoods of essential curves on $S$. Let $\varphi : D^2(S) \rightarrow D^2(S)$ be an automorphism of $D^2(S)$. Then there exists an automorphism $\tau : C(S) \rightarrow C(S)$ such that $\varphi \circ i = i \circ \tau$. \label{lem:cinv} \end{lem}

\begin{proof}
Let $a$ be a vertex of $C(S)$, $x = i(a)$, and $u = \varphi(x)$. Note that $x$ is an annular vertex of $D^2(S)$. Since $\varphi \in Aut(D^2(S))$, it follows by Corollary \ref{cor:D2annuli}, that $u$ is an annular vertex of $D^2(S)$. Hence, there exists a vertex $b$ of $C(S)$ such that $i(b) = u$. Since $i: C(S) \rightarrow D(S)$ is injective, such a vertex is unique. It follows that the correspondence $a \mapsto b$ yields a well-defined function $\tau: C_0(S) \rightarrow C_0(S)$ such that $\varphi(i(a)) = i(\tau(a))$ for every vertex $a$ of $C(S)$. Since curves on $S$ are disjoint if and only if they have disjoint regular neighborhoods, it follows that $\tau : C_0(S) \rightarrow C_0(S)$ extends to a simplicial map $\tau : C(S) \rightarrow C(S)$. Since $\varphi(i(a)) = i(\tau(a))$ for every vertex $a$ of $C(S)$, it follows that $\varphi \circ i = i \circ \tau : C(S) \rightarrow D^2(S)$. This shows that there exists a simplicial map $\tau: C(S) \rightarrow C(S)$ such that $\varphi \circ i = i \circ \tau : C(S) \rightarrow D^2(S)$. Likewise, there exists a simplicial map $\sigma: C(S) \rightarrow C(S)$ such that $\varphi^{-1} \circ i = i \circ \sigma : C(S) \rightarrow D^2(S)$. 
Since $i$ is injective, it follows that $\sigma$ is an inverse for $\tau$. Hence, $\tau : C(S) \rightarrow C(S)$ is an automorphism of $C(S)$.   
\end{proof}

\begin{thm} Suppose that $S$ is not a sphere with four holes, a torus with at most two holes, or a closed surface of genus two. Then the natural homomorphism $\rho: \Gamma^*(S) \rightarrow Aut(D^2(S))$ corresponding to the action of $\Gamma^*(S)$ on $D^2(S)$ is an isomorphism (i.e. every automorphism of $D^2(S)$ is induced by a homeomorphism $S \rightarrow S$ which is uniquely defined up to isotopy on $S$).
\label{thm:autD2} \end{thm}

\begin{proof} We begin by showing that $\rho$ is surjective. To this end, we let $\varphi \in Aut(D^2(S))$ and show that there exists a homeomorphism $H : S \rightarrow S$ such that $\varphi = H_* : D^2(S) \rightarrow D^2(S)$.

To simplify the exposition, we identify $C(S)$, via $i : C(S) \rightarrow D(S)$, with its image in $D^2(S)$ under $i : C(S) \rightarrow D(S)$. Since $S$ is not a torus with one hole, using this identification, we may restate Lemma \ref{lem:cinv} as saying that $\varphi$ restricts to an element $\tau$ of $Aut(C(S))$. 

Since $S$ is neither a sphere with at most four holes nor a torus with at most two holes, it follows from Theorem 1 of Ivanov \cite{ivanov1} and Theorem 1 of Korkmaz \cite{korkmaz} (see Luo \cite{luo} for a different proof) that there exists a homeomorphism $H : S \rightarrow S$ such that $\tau = H_* : C(S) \rightarrow C(S)$. Let $\psi = H^{-1}_* \circ \varphi : D^2(S) \rightarrow D^2(S)$. Note that $\psi$ fixes every vertex of $C(S)$. 

We shall now show that $\psi$ is equal to the identity map of $D^2(S)$. That is to say, we shall show that $\varphi = H_* : D^2(S) \rightarrow D^2(S)$.

Let $v \in D^2(S)$. Since $\psi$ is an automorphism of $D^2(S)$ preserving $C(S)$, $\psi(Ann(v)) = Ann(\psi(v))$. On the other hand, since $Ann(v)$ is a subcomplex of $C(S)$ and $\psi$ fixes each vertex of $C(S)$, $\psi(Ann(v)) = Ann(v)$. 
Hence, $Ann(\psi(v)) = Ann(v)$. Since $S$ is not a sphere with four holes, a torus with at most two holes, or a closed surface of genus two, it follows from Proposition \ref{prop:ann3} that $\psi(v) = v$. 
This proves that $\varphi = H_* : D^2(S) \rightarrow D^2(S)$ and, hence, the natural homomorphism $\rho: \Gamma^*(S) \rightarrow D^2(S)$ is surjective. 

It remains to show that $\rho: \Gamma^*(S) \rightarrow D^2(S)$ is injective. To this end, suppose that $h$ is an element of the kernel of $\rho$. Let  
$H : S \rightarrow S$ be a homeomorphism representing $h$. Since $h \in ker(\rho)$, $H$ induces the trivial automorphism of $D^2(S)$. 

Let $\alpha$ be an essential curve on $S$ and $X$ be a regular neighborhood of $\alpha$ on $S$. It follows that $[X] = H_*[X] = [H(X)]$ and, hence, $H(X)$ is isotopic to $X$ on $S$. This implies that $H(\alpha)$ is isotopic to $\alpha$ on $S$. Thus $H : S \rightarrow S$ preserves the isotopy class of every essential curve on $S$. In other words, $h$ is in the kernel of the action of $\Gamma^*(S)$ on $D^2(S)$. 

Since $S$ is not a sphere with at most three holes, it follows from Proposition \ref{prop:orpres} that $H : S \rightarrow S$ is 
orientation-preserving. This implies that $h$ is in the kernel of the action of $\Gamma(S)$ on $D^2(S)$. Since $S$ is not a sphere with at most four holes, a torus with at most two holes, or a closed surface of genus two, it follows from \cite{ivmcc}, Lemma 5.1 and Theorem 5.3, that $h$ is equal to the identity element of $\Gamma^*(S)$. 

This proves that $\rho: \Gamma^*(S) \rightarrow D^2(S)$ is injective, completing the proof.

\end{proof}

\begin{rem} Theorem \ref{thm:autD2} is vacuously true for spheres with at most three holes. That it fails for spheres with four holes, tori with at most two holes, and closed surfaces of genus two, follows from  detailed descriptions of $Aut(D^2(S))$ for these surfaces given in a sequel to this paper, in which we conduct a more detailed study of the complex of domains \cite{mccpap}.  
\label{rem:autD2} \end{rem}

\section{Recognizing biperipheral edges in $D(S)$} 

In this section, we shall characterize biperipheral edges of $D(S)$ .

Each biperipheral pair of pants on $S$ has a unique biperipheral boundary component. It follows that there is a natural map from the set of vertices of $D(S)$ corresponding to biperipheral pairs of pants on $S$ to the set of vertices of $D(S)$ corresponding to biperipheral curves on $S$. This map is a bijection if and only if $S$ is not a sphere with four holes. 

\begin{defn} Suppose that $X$ and $Y$ are domains on $S$ such that $X$ is a biperipheral pair of pants on $S$ and $Y$ is isotopic to a regular neighborhood of the unique essential boundary component of $X$ on $S$. Then we say that {\it $\{X,Y\}$ is a biperipheral pair of domains on $S$} and the edge $\{[X],[Y]\}$ of $D(S)$ is a {\it biperipheral edge of $D(S)$}.\label{defn:biperipheral edge} \end{defn}

Suppose that $\{X,Y\}$ is a biperipheral pair of domains on $S$. We may assume that $X$ is a biperiphal pair of pants on $S$. Since $Y$ is isotopic to a regular neighborhood of an essential boundary component of $X$ on $S$, $Y$ is isotopic to a domain $Y_1$ on $S$ which is disjoint from $X$. Since $X$ is not an annulus and $Y_1$ is an annulus, $X$ and $Y_1$ are not isotopic on $S$. It follows that $\{X,Y_1\}$ is a system of domains on $S$ and, hence, $\{[X],[Y]\} = \{[X],[Y_1]\}$ is, indeed, an edge of $D(S)$. Since $\{X,Y_1\}$ is both a biperipheral pair of domains on $S$ and a system of domains on $S$, we say that {\it $\{X,Y_1\}$ is a biperipheral system of domains on $S$}.

\begin{prop}[vertices with nested stars in $D(S)$] Let $x$ and $y$ be distinct vertices of $D(S)$. Then the following are equivalent: 

\begin{enumerate} 

\item \label{nestedstars1} $St(x,D(S)) \subset St(y,D(S))$.

\item \label{nestedstars2} There exist disjoint domains, $X$ and $Y$, on $S$ representing $x$ and $y$ which belong to one of the following cases:

\begin{enumerate} 

\item \label{nestedstarsa} $X$ is not an annulus and $Y$ is an annulus on $S$ which is joined to $X$ by exactly one annular codomain of $X \cup Y$ on $S$.  

\item \label{nestedstarsb} $X$ is not an annulus and $Y$ is an annulus on $S$ which is joined to $X$ by exactly two annular codomains of $X \cup Y$ on $S$.   

\item \label{nestedstarsc} $Y$ is a biperipheral pair of pants on $S$ which is joined to $X$ by exactly one annular codomain of $X \cup Y$ on $S$.  

\item \label{nestedstarsd} $Y$ is a monoperipheral pair of pants on $S$ which is joined to $X$ by exactly two annular codomains of $X \cup Y$ on $S$.  

\item \label{nestedstarse} $Y$ is a nonperipheral pair of pants on $S$ which is joined to $X$ by exactly three annular codomains of $X \cup Y$ on $S$.  

\end{enumerate}

\end{enumerate} 

\label{prop:nestedstars} \end{prop}

\begin{figure}[!hbp]
  \begin{center}
  \scalebox{0.60}{\includegraphics{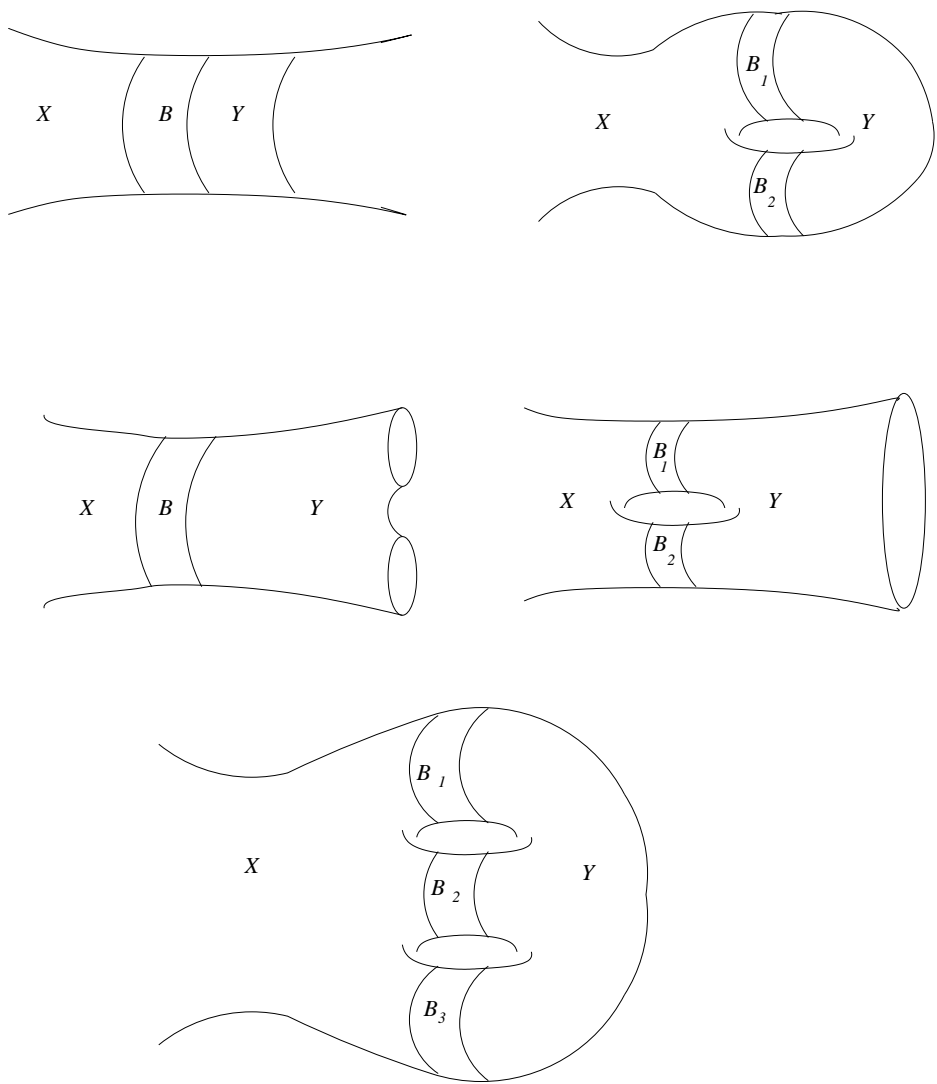}}
  \end{center}
   \caption{Ordered pairs of disjoint domains $(X,Y)$ representing ordered pairs of vertices $(x,y)$ that satisfy the equation $St(x,D(S))\subset St(y,D(S))$ with $x \neq y$. See Proposition \ref{prop:nestedstars}.}
  \label{fig:nestedstars}\end{figure}

\begin{proof} Suppose that $St(x,D(S)) \subset St(y,D(S))$. 

Since $x \in St(x,D(S))$, it follows that $x \in St(y,D(S))$. Since $x \neq y$, this implies that $\{x,y\}$ is an edge of $D(S)$. That is to say, $x$ and $y$ are represented by disjoint nonisotopic domains $X$ and $Y$ on $S$. 

Suppose that $Y$ is not elementary. By Proposition \ref{prop:nonelemdom}, 
there exist curves $\alpha$ and $\beta$ on $S$ such that $i(\alpha,\beta) \neq 0$ and $\alpha$ and $\beta$ are contained in the interior of $Y$. Let $W$ be a regular neighborhood of $\alpha$ on $S$ such that $W$ is contained in the interior of $Y$. Since $W$ is contained in $Y$ and $X$ is disjoint from $Y$, $X$ is disjoint from $W$. This implies that $\{x,w\}$ is a simplex of $D(S)$, where $w$ is the vertex of $D(S)$ represented by $W$. It follows that $w$ is a vertex of $St(x,D(S))$ and, hence, $w$ is a vertex of $St(y,D(S))$. That is to say, $\{y,w\}$ is a simplex of $D(S)$. Since $\{y,w\}$ is a simplex of $D(S)$, either $y = w$ or $\{y,w\}$ is an edge of $D(S)$. Since $Y$ is not an annulus on $S$ and $W$ is an annulus on $S$, $Y$ is not isotopic to $W$ on $S$. That is to say, $y \neq w$. Hence, $\{y,w\}$ is an edge of $D(S)$. It follows that $W$ is isotopic on $S$ to a domain on $S$ which is disjoint from $Y$. Since $\alpha$ is contained in $W$, it follows that $\alpha$ is isotopic on $S$ to a curve $\alpha_1$ which is disjoint from $Y$ and, hence, from $\beta$. Since $\alpha$ is isotopic on $S$ to $\alpha_1$ and $\alpha_1$ and $\beta$ are disjoint, it follows that $i(\alpha,\beta) = i(\alpha_1,\beta) = 0$, which is a contradiciton. Hence, $Y$ is elementary. 

Suppose that there exists an essential boundary component $\alpha$ of $Y$ which is not isotopic to any essential boundary component of $X$. 
Since $X$ and $Y$ are disjoint, it follows from Proposition \ref{prop:detectisotopy} that there exists an essential curve $\gamma$ on $S$ such that $i(\alpha,\gamma) \neq 0$ and $i(\beta,\gamma) = 0$ for every essential boundary component $\beta$ of $X$. 

We may assume that the collection $\mathcal{C}$ of curves on $S$ consisting of $\alpha$, $\gamma$, and the essential boundary components of $X$ on $S$ is in minimal position. It follows from the above constraints on geometric intersection numbers, that $\gamma$ is disjoint from $X$. 

Hence, there exists a regular neighborhood $Z$ of $\gamma$ on $S$ such that $Z$ is disjoint from $X$. Since $Z$ is disjoint from $X$, $Z$ represents a vertex $z$ of $St(x,D(S))$ and, hence, of $St(y,D(S))$. Thus, $\{y,z\}$ is a simplex of $D(S)$. 

Since $\{y,z\}$ is a simplex of $D(S)$, either $y = z$ or $\{y,z\}$ is an edge of $D(S)$. 

Suppose that $y = z$. That is to say, suppose that $Y$ is isotopic to $Z$ on $S$. Since $Z$ is an annulus on $S$, it follows that $Y$ is an annulus on $S$. Thus $Y$ is isotopic to a regular neighborhood of its essential boundary component $\alpha$. Since $Z$ is a regular neighborhood of $\gamma$ and $Y$ is isotopic to $Z$, it follows that $\alpha$ is isotopic to $\gamma$. Hence, $i(\alpha,\gamma) = i(\alpha,\alpha) = 0$ which is a contradiction. 

\end{proof}

\begin{prop}[vertices with the same star in $D(S)$] Let $x$ and $y$ be distinct vertices of $D(S)$. Then the following are equivalent: 

\begin{enumerate} 

\item \label{star1} $St(x,D(S)) = St(y,D(S))$.

\item \label{star2} There exist disjoint domains, $X$ and $Y$, on $S$ representing $x$ and $y$ which belong to one of the following cases:

\begin{enumerate} 

\item \label{stara} $X$ is a biperipheral pair of pants on $S$, $Y$ is an annulus on $S$, and $X \cup Y$ has exactly two codomains, exactly one of which is an annulus joining $X$ to $Y$. 

\item \label{starb} Case (\ref{stara}) with the roles of $X$ and $Y$ interchanged.

\item \label{starc} $S$ is a sphere with four holes, $X$ and $Y$ are biperipheral pairs of pants on $S$, and $X \cup Y$ has exactly one codomain, an annulus joining $X$ to $Y$. 

\item \label{stard} $S$ is a torus with two holes, $X$ and $Y$ are monoperipheral pairs of pants on $S$, and $X \cup Y$ has exactly two codomains, both of which are annuli joining $X$ to $Y$.

\item \label{stare} $S$ is a closed surface of genus two, $X$ and $Y$ are pairs of pants on $S$, and $X \cup Y$ has exactly three codomains, all of which are annuli joining $X$ to $Y$.

\item \label{starf} $S$ is a torus with one hole, $X$ is a monoperipheral pair of pants on $S$, $Y$ is an annulus on $S$, and $X \cup Y$ has exactly two codomains, both of which are annuli joining $X$ to $Y$.

\item \label{starg} Case (\ref{starf}) with the roles of $X$ and $Y$ interchanged.

\end{enumerate} 

\end{enumerate} 

\label{prop:star} \end{prop}

\begin{figure}[!hbp]
  \begin{center}
  \scalebox{0.60}{\includegraphics{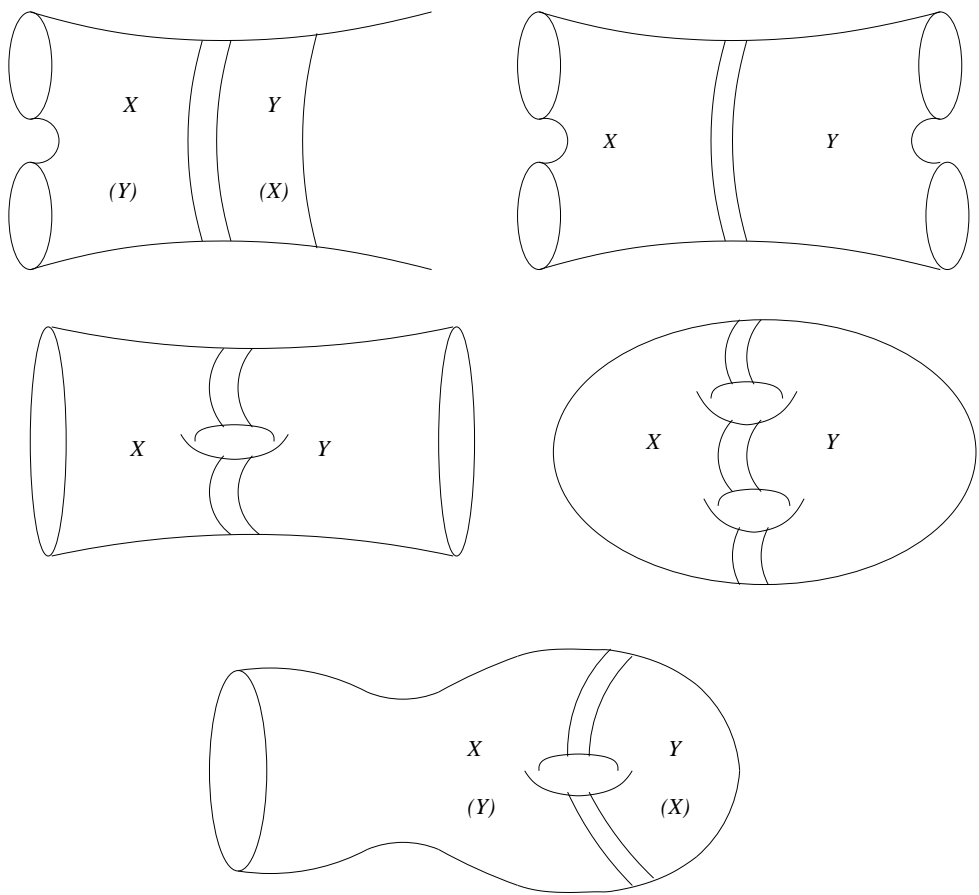}}
  \end{center}
   \caption{The seven topological types of ordered pairs of disjoint domains $(X,Y)$ representing ordered pairs of vertices $(x,y)$ that satisfy the equation $St(x,D(S))=St(y,D(S))$ with $x \neq y$. See Proposition \ref{prop:star}.}
  \label{fig:figstar}\end{figure}

\begin{proof} 

We begin by proving that (\ref{star1}) implies (\ref{star2}). To this end, suppose that $St(x,D(S)) = St(y,D(S))$. Since $x \neq y$ and $St(x,D(S)) \subset St(y,D(S))$, it follows from Proposition \ref{prop:nestedstars} that $x$ and $y$ are represented by disjoint nonisotopic domains $X$ and $Y$ satisfying one of the five cases, (\ref{nestedstarsa}), (\ref{nestedstarsb}), (\ref{nestedstarsc}), (\ref{nestedstarsd}), or  (\ref{nestedstarse}), of Proposition \ref{prop:nestedstars}. 

Note that for such domains, $X$ and $Y$, the ordered triple $(S,X,Y)$ is uniquely determined up to isotopies on $S$. Since $x \neq y$ and $St(y,D(S)) \subset St(x,D(S))$, it follows that $(Y,X)$ satisfies one of the five cases obtained by interchanging the roles of $X$ and $Y$ in Proposition \ref{prop:nestedstars}. 

Hence, $X$ and $Y$ are each either an annulus or a pair of pants. Moreover, if either is a pair of pants, all of its essential boundary components are joined to essential boundary components of the other by annular codomains of $X \cup Y$ and if either is an annulus, one or both of its essential boundary components is joined to the other by annular codomains of $X \cup Y$. Since $X$ is not isotopic to $Y$ on $S$ and $X$ and $Y$ are joined by at least one annular codomain of $X \cup Y$ on $S$, it follows that either $X$ is a pair of pants on $S$ or $Y$ is a pair of pants on $S$.

If $X$ and $Y$ are both pairs of pants, it follows that they have the same number $n$ of essential boundary components, with $1 \leq n \leq 3$. Hence, $X$ and $Y$ satisfy Case (\ref{starc}), when $n = 1$; Case (\ref{stard}), when $n = 2$; and Case (\ref{stare}), when $n = 3$. 

If $X$ is a pair of pants and $Y$ is an annulus, it follows that $X$ is either biperipheral when $X$ is joined to $Y$ by exactly one annular codomain of $X \cup Y$ on $S$ or monoperipheral when $X$ is joined to $Y$ by exactly two annular codomains of $X \cup Y$ on $S$. Hence, $X$ and $Y$ satisfy Case (\ref{stara}), when $X$ is biperipheral; and Case (\ref{starf}), when $X$ is monoperipheral. 

Likewise, if $X$ is an annulus and $Y$ is a pair of pants, then $X$ and $Y$ satisfy Case (\ref{starb}) or Case (\ref{starg}), according as $Y$ is either biperipheral or monoperipheral. 

This proves that (\ref{star1}) implies (\ref{star2}).

We now prove that (\ref{star2}) implies (\ref{star1}). 

To this end, suppose that $X$ and $Y$ are disjoint domains as in (\ref{star2}). We must prove that $Star(x,D(S)) = Star(y,D(S))$. Since $D(S)$ is a flag complex, it follows from Proposition \ref{prop:flagstar} that it suffices to prove that $St_0(x,D(S)) = St_0(y,D(S))$. 

We shall give the arguments for Cases (\ref{stara}) and (\ref{starc}). The argument for Case (\ref{starb}) is similar to that for Case (\ref{stara}). The argument for each of Cases (\ref{stard}), (\ref{stare}), (\ref{starf}), and (\ref{starg}) is similar to that for Case (\ref{starc}).

First, consider Case (\ref{stara}). Let $X$ and $Y$ be as in this case.

Suppose, on the one hand, that $w$ is a vertex of $St(x,D(S))$. In other words, suppose that $\{x,w\}$ is a simplex of $D(S)$. Then either $w = x$ or $\{x,w\}$ is an edge of $D(S)$. If $w = x$, then $\{y,w\}$ is the simplex $\{x,y\}$ of $D(S)$. Suppose that $\{x,w\}$ is an edge of $D(S)$. Then $w$ is represented by a domain $W$ on $S$ which is disjoint from $X$. Since $W$ is a domain on $S$ disjoint from $X$ and $Y$ is an annulus on $S$ which is joined to $X$ along the unique essential boundary component of $X$ by the unique annular codomain of $X \cup Y$ on $S$, it follows that $Y$ is isotopic on $S$ to a domain which is disjoint from $W$. It follows, in any case, that $\{y,w\}$ is a simplex of $D(S)$. That is to say, $w$ is a vertex of $St(y,D(S))$. 

Suppose, on the other hand, that $w$ is a vertex of $St(x,D(S))$. In other words, suppose that $\{y,w\}$ is a simplex of $D(S)$. Then either $w = y$ or $\{y,w\}$ is an edge of $D(S)$. If $w = y$, then $\{x,w\}$ is the simplex $\{x,y\}$ of $D(S)$. Suppose that $\{y,w\}$ is an edge of $D(S)$. Then $w$ is represented by a domain $W$ on $S$ which is disjoint from $Y$. Since $W$ is a domain on $S$ disjoint from $Y$ and $Y$ is an annulus on $S$ which is joined to $X$ along the unique essential boundary component of $X$ by the unique annular codomain of $X \cup Y$ on $S$, it follows that $W$ is isotopic on $S$ to a domain on $S$ which is contained in either $X$ or the complement of $X$. If $W$ is contained in the complement of $X$, then $\{x,w\}$ is a simplex of $D(S)$. Suppose that $W$ is contained in $X$. Since $W$ is a domain on $S$ which is contained in the biperipheral pair of pants $X$ on $S$, either $W$ is isotopic to $X$ on $S$ or $W$ is isotopic to a regular neighborhood of the unique essential boundary component of $X$ on $S$ and, hence, to $Y$. Hence, $\{x,w\}$ is equal to either the simplex $\{x\}$ of $D(S)$ or the simplex $\{x,y\}$ of $D(S)$. It follows, in any case, that $\{x,w\}$ is a simplex of $D(S)$. That is to say, $w$ is a vertex of $St(x,D(S))$. 

This proves that $St_0(x,D(S)) = St_0(y,D(S))$. This completes the argument for Case (\ref{stara}).

Now consider Case (\ref{starc}). Let $X$ and $Y$ be as in this case. Let $Z$ be the unique codomain of $X \cup Y$ on $S$. Note that $Z$ is an annulus on $S$ joining the unique essential boundary component of $X$ on $S$ to the unique essential boundary component of $Y$ on $S$. Hence, $S = X \cup Z \cup Y$. 

Since $X$ and $Y$ are biperipheral pairs of pants on $S$, it follows that $X \cup Z$ and $Y \cup Z$ are biperipheral pairs of pants on $S$ which are isotopic to $X$ and $Y$ on $S$ and $Z$ is isotopic on $S$ to regular neighborhoods on $S$ of the unique essential boundary components of each of $X$ and $Y$.

Suppose, on the one hand, that $w$ is a vertex of $St(x,D(S))$. In other words, suppose that $\{x,w\}$ is a simplex of $D(S)$. Then either $w = x$ or $\{x,w\}$ is an edge of $D(S)$. If $w = x$, then $\{y,w\}$ is the simplex $\{x,y\}$ of $D(S)$. Suppose that $\{x,w\}$ is an edge of $D(S)$. Then $w$ is represented by a domain $W$ on $S$ which is disjoint from $X$. Since $W$ is a domain on $S$ disjoint from $X$, it follows that $W$ is a domain on $S$ contained in $Y \cup Z$. Since $Y \cup Z$ is isotopic on $S$ to $Y$, we may assume that $W$ is contained in $Y$. Since $Y$ is a biperipheral pair of pants on $S$, it follows that $W$ is isotopic on $S$ to either $Y$ or the unique essential boundary component of $Y$ on $S$ and, hence, to $Z$. In other words, either $w = y$ or $w = z$. Hence, $\{y,w\}$ is equal to either the simplex $\{y\}$ of $D(S)$ or the simplex $\{y,x\}$ of $D(S)$. This shows, in any case, that $\{y,w\}$ is a simplex of $D(S)$. That is to say, $w$ is a vertex of $St(y,D(S))$. 

This proves that $St_0(x,D(S)) \subset St_0(y,D(S))$. By interchanging the roles of $X$ and $Y$, it follows that $St_0(y,D(S)) \subset St_0(x,D(S))$. Hence, $St_0(x,D(S)) = St_0(y,D(S))$. This completes the argument for Case (\ref{starc}). 

Since, as indicated above, the remaining cases follow by similar arguments, it follows that (\ref{star2}) implies (\ref{star1}).

\end{proof}

\begin{prop} Suppose that $S$ is not a sphere with four holes. Let $\{x,y\}$ be a pair of distinct vertices of $D(S)$. Let $\varphi \in Aut(D(S))$. Then $\{x,y\}$ is a biperipheral edge if and only if $\{\varphi(x),\varphi(y)\}$ is a biperipheral edge. 
\label{prop:biperipheraledges} \end{prop}

\begin{proof} Suppose, on the one hand, that $\{x,y\}$ is a biperipheral edge of $D(S)$.
(Note that since $\{x,y\}$ is a biperipheral edge of $D(S)$, $S$ has at least two boundary components.) 

We may assume that $x$ and $y$ are represented by disjoint domains $X$ and $Y$ on $S$ satisfying case (\ref{stara}) of Proposition \ref{prop:star}. Let $A$ be the unique codomain of $X \cup Y$ which joins $X$ to $Y$. Then $X \cup A \cup Y$ is a biperipheral pair of pants $P$ on $S$ which is isotopic on $S$ to $X$. Since $S$ is not a sphere with four holes, the unique codomain $Q$ of $P$ on $S$ is nonelementary. Since $P$ represents the vertex $x$ of $D(S)$, it follows that $Lk(x,D(S))$ has infinitely many vertices. 

 By Proposition \ref{prop:star}, $St(x,D(S)) = St(y,D(S))$. Since $\varphi : D(S) \rightarrow D(S)$ is an automorphism of $D(S)$, it follows that $\{\varphi(x),\varphi(y)\}$ is an edge of $D(S)$, $Lk(\varphi(x),D(S))$ has infinitely many vertices, and $St(\varphi(x),D(S)) = St(\varphi(y),D(S))$. 

Since $\varphi(x) \neq \varphi(y)$ and $St(\varphi(x),D(S)) = St(\varphi(y),D(S))$, it follows from Proposition \ref{prop:star} that $x$ and $y$ are represented by disjoint domains $X'$ and $Y'$ satisfying one of the seven cases of Proposition \ref{prop:star}. 

Suppose that $\{\varphi(x),\varphi(y)\}$ is not a biperipheral edge of $D(S)$. Then $X'$ and $Y'$ satisfy one of cases (\ref{starc}), (\ref{stard}), (\ref{stare}), (\ref{starf}), or (\ref{starg}) of Proposition \ref{prop:star}. Note that in any case, since $X'$ represents the vertex $\varphi(x)$ of $D(S)$, it follows that $Lk(\varphi(x), D(S))$ has at most four vertices, which is a contradiction. (In fact, $Lk(\varphi(x), D(S))$ has at most three vertices.) Hence, $\{\varphi(x),\varphi(y)\}$ is a biperipheral edge of $D(S)$. 

Suppose, on the other hand, that $\{\varphi(x),\varphi(y)\}$ is a biperipheral edge of $D(S)$. Then, since $\varphi^{-1} : D(S) \rightarrow D(S)$ is an automorphism of $D(S)$, it follows from the above argument, that $\{x,y\}$ is a biperipheral edge of $D(S)$. 

This completes the proof. 

\end{proof}

\section{Exchange automorphisms of $D(S)$} 

Throughout the rest of this paper, let $\mathcal{E}$ denote the set of biperipheral edges of $D(S)$. 

\begin{prop} Suppose that $S$ is not a sphere with four holes. Then there exists a monomorphism $\Phi : \mathcal{B}(\mathcal{E}) \rightarrow Aut(D(S))$ from the Boolean algebra $\mathcal{B}(\mathcal{E})$ of all subsets of $\mathcal{E}$ to $Aut(D(S))$ such that for each collection $\mathcal{F}$ of biperipheral edges of $D(S)$, $\Phi(\mathcal{F}) = \varphi_{\mathcal{F}}$ exchanges the two vertices of each biperipheral edge in $\mathcal{F}$ and fixes every vertex of $D(S)$ which is not a vertex of some biperipheral edge in $\mathcal{F}$.  
\label{prop:Dboolean} \end{prop}

\begin{proof} It follows from Propositions \ref{prop:star} and \ref{prop:xyanedge} that $\mathcal{E}$ is a collection of exchangeable edges of $D(S)$. Since $S$ is not a sphere with four holes, no two distinct edges in $\mathcal{E}$ have a common vertex. Hence, the result follows from Proposition \ref{prop:boolean}. 
\end{proof}

Following the language of Definition \ref{defn:boolean}, we call the image of the Boolean algebra $\mathcal{B}(\mathcal{E})$ under the monomorphism $\Phi$ of Proposition \ref{prop:Dboolean} the {\it Boolean subgroup $B_{\mathcal{E}}$ of $D(S)$}. In particular, the Boolean subgroup $B_{\mathcal{E}}$ is naturally isomorphic to the Boolean algebra $\mathcal{B}(\mathcal{E})$.

\begin{prop} Let $\varphi \in Aut(D(S))$, $\mathcal{F} \subset \mathcal{E}$ and $\mathcal{G} = \varphi(\mathcal{F})$. Then $\mathcal{G} \subset \mathcal{E}$ and $\varphi \circ \Phi_{\mathcal{F}} \circ \varphi^{-1} = \Phi_{\mathcal{G}}$.
\label{prop:DBooleanconjugation} \end{prop} 

\begin{proof} This is an immediate consequence of Propositions \ref{prop:biperipheraledges} and \ref{prop:booleanconjugation}.
\end{proof}

\begin{prop} $B_{\mathcal{E}}$ is a normal subgroup of $Aut(D(S))$. 
  \label{prop:DBooleannormal} \end{prop} 
  
  \begin{proof} This is an immediate consequence of Proposition \ref{prop:DBooleanconjugation}.
\end{proof}

\begin{prop} The monomorphism $\Phi : \mathcal{B}(\mathcal{E}) \rightarrow Aut(D(S))$ is natural with respect to the action of the extended mapping class group $\Gamma^*(S)$ on $D(S)$. More precisely, if $h \in \Gamma^*(S)$ and $\mathcal{F} \subset \mathcal{E}$, then $\Phi(h_*(\mathcal{F})) = h_* \circ \Phi(\mathcal{F}) \circ h_*^{-1}$. 
 \label{prop:DBooleanequivariance} \end{prop} 

\begin{proof} This is an immediate consequence of Propositions \ref{prop:DBooleanconjugation} and \ref{prop:DBooleannormal}.
\end{proof}

\begin{prop} There is a natural monomorphism:

$$\rho: \mathcal{B}(\mathcal{E}) \rtimes \Gamma^*(S) \longrightarrow Aut(D(S))$$
    
\noindent corresponding to the action of $\Gamma^*(S)$ on $D(S)$ and the induced action on the set $\mathcal{E}$ of biperipheral edges of $D(S)$.
\label{prop:semdirprod} \end{prop}

\begin{proof} Since, by Proposition \ref{prop:DBooleanequivariance}, the monomorphism $\Phi : \mathcal{B}(\mathcal{E}) \rightarrow B_{\mathcal{E}}$ is natural, there exists a natural homomorphism $\rho : \mathcal{B}(\mathcal{E}) \rtimes \Gamma^*(S) \longrightarrow Aut(D(S))$. 
Since a pair of pants is not homeomorphic to an annulus, a geometric automorphism of $D(S)$ cannot exchange the vertices of any biperipheral edge of $D(S)$. It follows that the image of the extended mapping class group $\Gamma^*(S)$ in $Aut(D(S))$ under the natural homomorphism $\rho : \Gamma^*(S) \rightarrow Aut(D(S))$ corresponding to the action of $\Gamma^*(S)$ on $D(S)$ has trivial intersection with the Boolean subgroup $B_{\mathcal{E}}$ of $Aut(D(S))$. 
Since the natural homomorphism $\Phi : \mathcal{B}(\mathcal{E}) \rightarrow B_{\mathcal{E}}$ is injective, it remains only to show that $\rho: \Gamma^*(S) \rightarrow Aut(D(S))$ is injective. To this end, suppose that $h \in \Gamma^*(S)$ is in the kernel of $\rho$. Since $D^2(S)$ is a subcomplex of $D(S)$, it follows that $h$ induces the trivial automorphism of $D^2(S)$. Since $S$ is not a sphere with four holes, a torus with at most two holes, or a closed surface of genus two, it follows from Theorem \ref{thm:autD2} that $h$ is equal to the identity element of $\Gamma^*(S)$. This proves that $\rho: \Gamma^*(S) \rightarrow Aut(D(S))$ is injective, completing the proof.
\end{proof}

\section{Automorphisms of $D(S)$}

Throughout this section, let $\mathcal{E}$ denote the set of biperipheral edges of $D(S)$. 

\begin{prop} Suppose that $S$ is not a sphere with four holes. Let $\pi: D(S) \rightarrow D^2(S)$ be the natural projection from $D(S)$ to $D^2(S)$ sending each vertex of $D(S)$ corresponding to a biperipheral pair of pants on $S$ to the annular vertex of $D(S)$ corresponding to its unique essential boundary component on $S$. If $\varphi \in Aut(D(S))$, then there exists a unique simplicial automorphism $\varphi_* : D^2(S) \rightarrow D^2(S)$ such that $\varphi_* \circ \pi = \pi \circ \varphi : D(S) \rightarrow D^2(S)$ .  
\label{prop:Dphi-star} \end{prop}

\begin{proof} Let $i : D^2(S) \rightarrow D^2(S)$ denote the inclusion map of the subcomplex $D^2(S)$ of $D(S)$ into $D(S)$ and $\varphi_* = \pi \circ \varphi \circ i : D^2(S) \rightarrow D^2(S)$. Note that $\varphi_* : D^2(S) \rightarrow D^2(S)$ is a simplicial map from $D^2(S)$ to $D^2(S)$. 

We shall prove that $\varphi_* \circ \pi = \pi \circ \varphi : D(S) \rightarrow D^2(S)$. To this end, let $x$ be a vertex of $D(S)$. 

Suppose, on the one hand, that $x \in D^2(S)$. Then, by the definition of $\pi : D(S) \rightarrow D^2(S)$, $\pi(x) = x$ and, hence, $(\varphi_* \circ \pi)(x) = \varphi_*(\pi(x)) = \varphi_*(x) = (\pi \circ \varphi \circ i)(x) = \pi(\varphi(i(x))) = \pi(\varphi(x)) = (\pi \circ \varphi)(x)$.

Suppose, on the other hand, that $x$ is not in $D^2(S)$. Since $x \in D(S)$ and $x$ is not in $D^2(S)$, $x$ is represented by a biperipheral pair of pants $X$ on $S$. Let $Y$ be a regular neighborhood of the unique essential boundary component of $X$ on $S$ and $y$ be the vertex of $D(S)$ represented by $Y$. Then $\{x,y\}$ is a biperipheral edge of $D(S)$. It follows from the definition of $\pi : D(S) \rightarrow D(S)$, that $\pi(x) = y = \pi(y)$. Moreover, it follows from Proposition \ref{prop:biperipheraledges} that $\{\varphi(x),\varphi(y)\}$ is a biperipheral edge of $D(S)$. 
Hence, either $\varphi(x)$ is represented by a biperipheral pair of pants on $S$ or $\varphi(y)$ is represented by a biperipheral pair of pants on $S$. In the former case,it follows from the definition of $\pi: D(S) \rightarrow D^2(S)$, that $\pi(\varphi(x)) = \varphi(y) = \pi(\varphi(y))$.  In the latter case,it follows from the definition of $\pi: D(S) \rightarrow D^2(S)$, that $\pi(\varphi(x)) = \varphi(x) = \pi(\varphi(y))$. Hence, in any case, $\pi(\varphi(x)) = \pi(\varphi(y))$. It follows that $(\varphi_* \circ \pi)(x) = \varphi_*(\pi(x)) = \pi(\varphi(i(\pi(x))) = \pi(\varphi(\pi(x)) = 
\pi(\varphi(y)) = \pi(\varphi(x)) = (\pi \circ \varphi)(x)$. 

This shows, in any case, that $(\varphi_* \circ \pi)(x) = (\pi \circ \varphi)(x)$ and, hence, $\varphi_* \circ \pi = \pi \circ \varphi : D(S) \rightarrow D^2(S)$. 

Suppose that $\beta : D^2(S) \rightarrow D^2(S)$ is a simplicial map such that $\beta \circ \pi = \pi \circ \varphi : D(S) \rightarrow D^2(S)$. Then $\beta \circ \pi = \varphi_* \circ \pi : D(S) \rightarrow D^2(S)$. Since $\pi : D(S) \rightarrow D^2(S)$ is surjective, it follows that $\beta = \varphi_* : D(S) \rightarrow D^2(S)$. This proves that there exists a unique simplicial map $\varphi_* : D^2(S) \rightarrow D^2(S)$ such that $\varphi_* \circ \pi = \pi \circ \varphi : D(S) \rightarrow D^2(S)$. 

It remains only to prove that $\varphi_* : D^2(S) \rightarrow D^2(S)$ is a simplicial automorphism of $D^2(S)$. To this end, consider the simplicial automorphism $\psi = \varphi^{-1}: D(S) \rightarrow D(S)$. of $D(S)$. Repeating the above argument, we conclude that there exists a unique simplicial map $\psi_* : D^2(S) \rightarrow D^2(S)$ such that $\psi_* \circ \pi = \pi \circ \psi : D(S) \rightarrow D^2(S)$. 

It follows that $(\varphi_* \circ \psi_*) \circ \pi = \varphi_* \circ (\psi_* \circ \pi) = \varphi_* \circ (\pi \circ \psi) = (\varphi_* \circ \pi) \circ \psi = (\pi \circ \varphi) \circ \psi = \pi \circ (\varphi \circ \psi) = \pi : D(S) \rightarrow D^2(S)$. Since $\pi : D(S) \rightarrow D^2(S)$ is surjective, it follows that $\varphi_* \circ \psi_* : D^2(S) \rightarrow D^2(S)$ is the identity map of $D^2(S)$. Likewise, we conclude that $\psi_* \circ \varphi_* : D^2(S) \rightarrow D^2(S)$ is the identity map of $D^2(S)$. Hence, $\varphi_* : D^2(S) \rightarrow D^2(S)$ and $\psi_* : D^2(S) \rightarrow D^2(S)$ are inverse simplicial maps. This shows that $\varphi_* : D^2(S) \rightarrow D^2(S)$ is a simplicial automorphism of $D^2(S)$, completing the proof.  

\end{proof}

\begin{prop} Suppose that $S$ is not a sphere with four holes. Let $\pi: D(S) \rightarrow D^2(S)$ be the natural projection from $D(S)$ to $D^2(S)$ sending each vertex of $D(S)$ corresponding to a biperipheral pair of pants on $S$ to the annular vertex of $D(S)$ corresponding to its unique essential boundary component on $S$. Then there exists a unique homomorphism $\rho : Aut(D(S)) \rightarrow Aut(D^2(S))$ such that for each automorphism $\varphi \in Aut(D(S))$, $\rho(\varphi)$ is the unique simplicial automorphism $\varphi_* : D^2(S) \rightarrow D^2(S)$ such that $\varphi_* \circ \pi = \pi \circ \varphi : D(S) \rightarrow D^2(S)$. Moreover, there exists a natural exact sequence: 

$$1 \longrightarrow  B_{\mathcal{E}}  \longrightarrow 
    Aut(D(S))  \longrightarrow   Aut(D^2(S)) 
    \label{diagram8}.$$
    
\label{prop:Dexactsequence} \end{prop}

\begin{proof} By Proposition \ref{prop:Dphi-star}, there is a map $\rho :  Aut(D(S)) \rightarrow Aut(D^2(S))$ such that for each automorphism $\varphi \in Aut(D(S))$, $\rho(\varphi)$ is the unique simplicial automorphism $\varphi_* : D^2(S) \rightarrow D^2(S)$ such that $\varphi_* \circ \pi = \pi \circ \varphi : D(S) \rightarrow D^2(S)$. Suppose that $\varphi : D(S) \rightarrow D(S)$ and $\psi: D(S) \rightarrow D(S)$ are elements of $Aut(D(S))$. Since $\varphi_* : D^2(S) \rightarrow D^2(S)$ and $\psi_* : D^2(S) \rightarrow D^2(S)$ are automorphisms of $D^2(S)$,  $(\varphi_* \circ \psi_*) \circ \pi = \varphi_* \circ (\psi_* \circ \pi) = \varphi_* \circ (\pi \circ \psi) = (\varphi_* \circ \pi) \circ \psi = (\pi \circ \varphi) \circ \psi = \pi \circ (\varphi \circ \psi)$. It follows from the uniqueness clause of Proposition \ref{prop:Dphi-star} that $(varphi \circ \psi)_* = \varphi_* \circ \psi_* : D^2(S) \rightarrow D^2(S)$ and, hence, $\rho :  Aut(D(S)) \rightarrow Aut(D^2(S))$ is a homomorphism. This proves the existence and uniqueness of such a homomorphism $\rho :  Aut(D(S)) \rightarrow Aut(D^2(S))$.

Since $B_{\mathcal{E}}$ is by definition a subgroup of $Aut(D(S))$, the natural homomorphism $B_{\mathcal{E}} \rightarrow Aut(D(S))$ is injective. 

Suppose that $\varphi \in B_{\mathcal{E}}$. By the definition of $B_{\mathcal{E}}$, there exists a unique subset $\mathcal{F}$ of the collection $\mathcal{E}$ such that $\varphi$ exchanges the two vertices of each pair of distinct vertices of $D(S)$ in the collection $\mathcal{F}$ and fixes every vertex of $D(S)$ which is not one of the two vertices of some pair of distinct vertices of $D(S)$ in the collection $\mathcal{F}$. 

Suppose that $z$ is a vertex of $D^2(S)$. Since $\pi: D(S) \rightarrow D^2(S)$ is a surjective simplicial map, there exists a vertex $x$ of $D(S)$ such that $\pi(x) = z$.

Suppose, on the one hand, that $x$ is one of the two vertices of some distinct pair of vertices $\{x,y\}$ of $D(S)$ in the collection $\mathcal{F}$. Since $\{x,y\}$ is in $\mathcal{F}$, $\varphi$ interchanges $x$ and $y$ and, hence, $\varphi(x) = y$. Since $\mathcal{F}$ is a subset of $\mathcal{E}$, it follows from the definition of $\pi : D(S) \rightarrow D^2(S)$ that $\pi(y) = \pi(x)$. Hence, $\varphi_*(z) = \varphi_*(\pi(x)) = \pi(\varphi(x)) = 
\pi(y) = \pi(x) = z$.   

Suppose, on the other hand, that $x$ is not one of the two vertices of any distinct pair of vertices of $D(S)$ in the collection $\mathcal{F}$. Then, $\varphi(x) = x$. Hence, $\varphi_*(z) = \varphi_*(\pi(x)) = \pi(\varphi(x)) = \pi(x) = z$.   

In any case, it follows that the simplicial automorphism $\varphi_* : D^2(S) \rightarrow D^2(S)$ of $D^2(S)$ fixes every vertex of $D^2(S)$ and is, hence, the identity map of $D^2(S)$. This proves that the image of the natural homomorphism $B_{\mathcal{E}} \rightarrow Aut(D(S))$ is in the kernel of $\rho :  Aut(D(S)) \rightarrow Aut(D^2(S))$. 

Conversely, suppose that $\varphi : D(S) \rightarrow D(S)$ is in the kernel of  $\rho : Aut(D(S)) \rightarrow Aut(D^2(S))$.
Then $\varphi_* : D^2(S) \rightarrow D^2(S)$ is the identity map of $D^2(S)$. Let $x$ be a vertex of $D(S)$, $y = \varphi(x)$, and $z = \pi(x)$. Since $z$ is a vertex of $D^2(S)$, it follows that $z = \varphi_*(z) = \varphi_*(\pi(x)) = \pi(\varphi(x)) = \pi(y)$. 

Hence, $x$ and $y$ are in the same fiber of $\pi : D(S) \rightarrow D^2(S)$. It follows from the definition of $\pi : D(S) \rightarrow D^2(S)$ that either $y = x$ or $\{x,y\}$ is a pair of distinct vertices of $D(S)$ in the collection $\mathcal{E}$. 

Suppose that $y$ is not equal to $x$. Then $\{x,y\}$ is a pair of distinct vertices of $D(S)$ in the collection $\mathcal{E}$. Let $w = \varphi(y)$. Repeating the previous argument, with $(y,w,z)$ rather than $(x,y,z)$. we conclude that either $w = y$ or $\{y,w\}$ is a pair of distinct vertices of $D(S)$ in the collection $\mathcal{E}$. Suppose that $w = y$. Then $\varphi(y) = w = y = \varphi(x)$. Since $\varphi : D(S) \rightarrow D(S)$ is a simplicial automorphism and, hence, injective, it follows that $y = x$, which is a contradiction. It follows that $w$ is not equal to $y$ and, hence, $\{y,w\}$ is a pair of distinct vertices of $D(S)$ in the collection $\mathcal{E}$. Since  $\{x,y\}$ and $\{y,w\}$ are both pairs of distinct vertices of $D(S)$ in the collection $\mathcal{E}$ having at least one common vertex $y$ and no two distinct pairs of vertices in the collection $\mathcal{E}$ have a common vertex, it follows that $\{x,y\} = \{y,w\}$. Since $w$ is not equal to $y$, it follows that $w = x$. That is to say, $\varphi(y) = x$. 

This proves that for each vertex $x$ of $D(S)$, either $\varphi(x) = x$ or $x$ is one of the two vertices of a pair $\{x,y\}$ of distinct vertices of $D(S)$ in $\mathcal{E}$ and $\varphi$ exchanges $x$ and $y$.  

Let $\mathcal{F}$ be the subset of $\mathcal{E}$ consisting of all pairs of distinct vertices of $D(S)$ in $\mathcal{E}$ which are exchanged by $\varphi$. It follows that $\varphi : D(S) \rightarrow D(S)$ is equal to the generalized exchange $\varphi_{\mathcal{F}} : D(S) \rightarrow D(S)$ of $D(S)$ associated to $\mathcal{F}$. By the definition of $B_{\mathcal{E}}$, $\varphi$ is an element of $B_{\mathcal{E}}$. Hence, the kernel of $\rho : Aut(D(S)) \rightarrow Aut(D^2(S))$ is in the image of the natural homomorphism $B_{\mathcal{E}} \rightarrow Aut(D(S))$. This proves that the image of the natural homomorphism $B_{\mathcal{E}} \rightarrow Aut(D(S))$ is equal to the kernel of $\rho : Aut(D(S)) \rightarrow Aut(D^2(S))$. 
  
\end{proof}

\begin{thm} Suppose that $S$ is not a sphere with at most four holes, a torus with at most two holes, or a closed surface of genus two. Every automorphism of $D(S)$ is a composition of an exchange automorphism of $D(S)$ with a geometric automorphism of $D(S)$. 

More precisely, let $\mathcal{E}$ be the collection of biperipheral edges of $D(S)$. Let $\varphi \in Aut(D(S))$. Then there exists a unique subset $\mathcal{F}$ of $\mathcal{E}$ and a unique element $h$ of $\Gamma^*(S)$ such that $\varphi$ is equal to the composition $\varphi_{\mathcal{F}} \circ h_*$ of the exchange automorphism $\varphi_{\mathcal{F}}$ of $D(S)$ corresponding to $\mathcal{F}$ and the geometric automorphism $h_*$ of $D(S)$ induced by $h$. 
\label{thm:exchangegeom} \end{thm}

\begin{proof} Let $\varphi \in Aut(D(S))$. We begin by proving the existence of such a factorization of $\varphi$. Since $S$ is not a sphere with four holes,  it follows from Proposition \ref{prop:Dphi-star} that there exists a unique simplicial automorphism $\psi : D^2(S) \rightarrow D^2(S)$ such that $\psi \circ \pi = \pi \circ \varphi : D(S) \rightarrow D^2(S)$. 

Since $S$ is not a sphere with at most four holes, a torus with at most two holes, or a closed surface of genus two, it follows from Theorem \ref{thm:autD2} that there exists a homeomorphism $H : S \rightarrow S$ such that $\psi([X]) = [H(X)]$ for every domain $X$ on $S$ which is not a biperipheral pair of pants.

Let $G = H^{-1} : S \rightarrow S$ and $G_* : D(S) \rightarrow D(S)$ be the geometric automorphism of $D(S)$ defined by the rule $G_*([X]) = [G(X)]$ for every domain $X$ on $S$.

Note that $G_* \circ \varphi : D(S) \rightarrow D(S)$ is an automorphism of $D(S)$. We shall now show that $G_* \circ \varphi$ is an exchange automorphism. 

Since $S$ is not a sphere with four holes and $\varphi \in Aut(D(S))$, it follows from Proposition \ref{prop:biperipheraledges} that $\varphi(\mathcal{E}) = \mathcal{E})$. 

Suppose that $X$ is a domain on $S$ which is not a biperipheral pair of pants or a regular neighborhood of a biperipheral curve. Note that $[X]$ is not a vertex of an edge in $\mathcal{E}$. Since 
$\varphi$ is an automorphism of $D(S)$ and $\varphi(\mathcal{E}) = \mathcal{E}$, it follows that $\varphi([X])$ is not a vertex of an edge in $\mathcal{E}$. Hence, $\varphi([X]) = [Y]$ where $Y$ is a domain on $S$ which is not a biperipheral pair of pants or a regular neighborhood of a biperipheral curve.

By the definition of the natural projection $\pi : D(S) \rightarrow D^2(S)$, $\pi([X]) = [X]$ and 
$\pi([Y]) = [Y]$. Hence, $\pi(\varphi([X]) = \varphi([X])$.

Hence, $\varphi([X]) = \pi(\varphi([X])) = \psi(\pi([X]) = \psi([X]) = [H(X)]$. 

It follows that $(\varphi \circ G_*)([X]) = \varphi[G(X)] = [H(G(X))] = [X]$. 

This shows that $\varphi \circ G_*$ fixes every vertex of $D(S)$ which is not a vertex of an edge in $\mathcal{E}$. 

By a similar argument, it follows that if $X$ is a domain on $S$ which is a biperipheral pair of pants and $Y$ is a regular neighborhood of the corresponding biperipheral curve, then either
(i) $(\varphi \circ G_*)([X]) = [X]$  and $(\varphi \circ G_*)([Y]) = [Y]$ or (ii) $(\varphi \circ G_*)([X]) = [Y]$ and $(\varphi \circ G_*)([Y]) = [X]$. 

Let $i = \varphi \circ G_*$. It follows that $i$ is an exchange automorphism of $D(S)$. 

Since $i = \varphi \circ G_*$, we conclude that $\varphi = \varphi_{\mathcal{F}} \circ h_*$ where $\mathcal{F}$ is the subcollection of $\mathcal{E}$ consisting of all biperipheral edges of $D(S)$ whose vertices are exchanged by $i$, and $h$ is the isotopy class of $H : S \rightarrow S$. 

This proves the existence of such a factorization of $\varphi$. 

Suppose that $\Phi_{\mathcal{F}} \circ h_* = \Phi_{\mathcal{P}} \circ q_*$. 
Then $\Phi_{\mathcal{P} \triangle \mathcal{F}} = (g \cdot h^{-1})_*$. 

Since an automorphism of $D(S)$ which is induced by a homeomorphism of $S$ cannot exchange an annular vertex with a nonannular vertex of $D(S)$, it follows that $\mathcal{P} \triangle \mathcal{F} = \emptyset$. In other words, $\mathcal{F} = \mathcal{P}$. 

Since $\Phi_{\mathcal{P} \triangle \mathcal{F}} = (g \cdot h^{-1})_*$ and 
$\mathcal{F} = \mathcal{P}$, it follows that $(g \cdot h^{-1})_*$ is the trivial automorphism $id : D(S) \rightarrow D(S)$ of $D(S)$. In other words, $g \cdot h^{-1}$ acts trivially on $D(S)$. 

Since $D^2(S)$ is a subcomplex of $D(S)$, it follows that $g \cdot h^{-1}$ acts trivially on $D^2(S)$. Since $S$ is not a sphere with four holes, a torus with at most two holes, or a closed surface of genus two, it follows from Theorem \ref{thm:autD2} that $g \cdot h^{-1}$ is equal to the identity element of $\Gamma^*(S)$. In other words, $g$ is equal to $h$. 

This proves the uniqueness of such a factorization of $\varphi$, completing the proof.

\end{proof}

We can summarize the preceding results as follows. 

\begin{thm} Suppose that $S$ is not a sphere with at most four holes, a torus with at most two holes, or a closed surface of genus two. Then we have a natural commutative diagram of exact sequences: 
 
\[ \begin{array}{cccccccccc}
    1 &\longrightarrow & \mathcal{B}(\mathcal{E})   & \longrightarrow &
    \mathcal{B}(\mathcal{E}) \rtimes \Gamma^*(S) & \longrightarrow  & \Gamma^*(S) & \longrightarrow & 1\\
   \   & \  & \simeq\big\downarrow & \  & \simeq\big\downarrow & \  & \simeq\big\downarrow
      & \ & \  \\
     1 & \longrightarrow & B_{\mathcal{E}}    & \longrightarrow &
    Aut(D(S)) & \longrightarrow  & Aut(D^2(S)) & \longrightarrow & 1 
    \label{diagram7} \end{array}\]

\label{thm:exchangegeom2} \end{thm}

\begin{proof} The exactness of the first row of the above diagram follows immediately from the definition of a semi-direct product. 

The commutativity of the left hand square follows from Propositions \ref{prop:Dboolean} and \ref{prop:semdirprod}. 

The commutativity of the right hand square follows from Propositions \ref{prop:semdirprod} and \ref{prop:Dphi-star}. 

The isomorphism $\mathcal{B}(\mathcal{E}) \stackrel{\simeq}{\rightarrow} B_{\mathcal{E}}$ follows from Proposition \ref{prop:Dboolean} and the definition of the Boolean subgroup $B_{\mathcal{E}}$ of $Aut(D(S))$. 

Since $S$ is not a sphere with at most four holes, a torus with at most two holes, or a closed surface of genus two, it follows from Theorem \ref{thm:autD2} that the natural homomorphism $\Gamma^*(S) \stackrel{\simeq}{\rightarrow} Aut(D^2(S))$ is an isomorphism.

Since the natural homomorphisms  $\mathcal{B}(\mathcal{E}) \rtimes \Gamma^*(S) \longrightarrow  \Gamma^*(S)$ and $\Gamma^*(S) \rightarrow Aut(D^2(S))$ are both surjective, it follows from the commutativity of the right hand square that the natural homomorphisms $Aut(D(S)) \rightarrow Aut(D^2(S))$ is also surjective. Hence, since $S$ is not a sphere with four holes, the exactness of the second row of the above diagram follows from Proposition \ref{prop:Dexactsequence}.

This shows that the diagram is a commutative diagram of exact sequences. Since the vertical arrows on the left and right are both isomorphisms, it follows from standard results that the natural monomorphism $\mathcal{B}(\mathcal{E}) \rtimes \Gamma^*(S) \rightarrow Aut(D(S))$ of Proposition \ref{prop:semdirprod} is an isomorphism, completing the proof.

\end{proof}


\begin{thebibliography}{XXXX}

\bibitem{behrstockmargalit} J. Behrstock and D. Margalit, Curve complexes and finite index subgroups of mapping class groups, LANL ArXiv, math.GT/0504328, to appear in Geometriae Dedicata

\bibitem{brendlemargalit} T. E. Brendle and D. Margalit, Commensurations of the Johnson kernel, Geometry \& Topology, vol. 8 (2004) no. 37,  1361-1384.

\bibitem{dehn} M. Dehn, Die Gruppe der Abbildungsklassen, Acta Math. {\bf 69}, 135-206 (1938)

\bibitem{farbivanov} B. Farb and N. V. Ivanov, The Torelli geometry and its applications: research announcement, Math. Res. Lett. 12 (2005), no. 2-3, 293-301.

\bibitem{flp} A. Fathi, F. Laudenbach. and V. Po\'{e}naru, {\it Travaux de Thurston  sur les surfaces}, S\'eminaire Orsay, Ast\'{e}risque, Vol. 66-67, Soc. Math. de France, 1979.

\bibitem{harvey} W. J. Harvey, Geometric structure of surface mapping class groups, {\it Homological group theory (Proc. Sympos., Durham, 1977)}, pp. 255--269, London Math. Soc. Lecture Note Ser., 36, Cambridge Univ. Press, Cambridge-New York, 1979. 

\bibitem{irmak1} E. Irmak, Superinjective simplicial maps of complexes of curves and injective homomorphisms of subgroups of mapping class groups, Topology 43 (2004), no. 3, 513-541.

\bibitem{irmak2} E. Irmak, Complexes of nonseparating curves and mapping class groups, LANL ArXiv, math.GT/0407285.

\bibitem{irmak3} E. Irmak, Superinjective simplicial maps of complexes of curves and injective homomorphisms of subgroups of mapping class groups. II, Topology Appl. 153 (2006), no. 8, 1309--1340.

\bibitem{irmakkorkmaz} E. Irmak and M. Korkmaz, Automorphisms of the Hatcher-Thurston complex, LANL ArXiv, math.GT/0409033, to appear in
Israel Journal of Mathematics 

\bibitem{ivanov1} N. V. Ivanov, Automorphisms of complexes of curves and of Teichm\"{u}ller spaces, Preprint IHES/M/89/60, 1989, 13 pp.; Also in: {\it Progress in knot theory and related topics,} Travaux en Cours, V. 56, {\it Hermann, Paris,} 1997, 113-120.

\bibitem{ivanov2} N. V. Ivanov, {\it Subgroups of Teichm\"{u}ller Modular Groups},  Translations of Mathematical Monographs, Vol. 115, American Math. Soc., Providence, Rhode Island,  1992.

\bibitem{ivmcc} N. V. Ivanov and J. D. McCarthy, On injective homomorphisms between Teichm\"{u}ller modular groups, I,  Invent. math. 135, 425-486 (1999).

\bibitem{johnson} D. L. Johnson, Homeomorphisms of a surface which act trivially on homology, Proc. Amer. Math. Soc. {\bf 75}, 119-125 (1979)

\bibitem{korkmaz} M. Korkmaz, Automorphisms of complexes of curves on punctured spheres and on punctured tori, Topology and its Applications  95 no. 2 (1999), 85-111.

\bibitem{luo} F. Luo,  Automorphisms of the complex of curves,
Topology 39 (2000), no. 2, 283-298.

\bibitem{margalit} D. Margalit, Automorphisms of the pants complex.
Duke Math. J. 121 (2004), no. 3, 457-479.

\bibitem{mm1} H. A. Masur and Y. N. Minsky, Geometry of the complex of curves I: Hyperbolicity, Invent. math. 138, 103-149 (1999) 

\bibitem{mcc3} J. D. McCarthy, Automorphisms of surface mapping class groups. A recent theorem of N. Ivanov, Invent. Math. V. 84, F. 1 (1986), 49-71.

\bibitem{mcc5} J. D. McCarthy, Simplicial representations of the mapping class group, to appear in the {\it Handbook of Teichm\"{u}ller Theory}, vol. II, European Mathematical Society

\bibitem{mccpap} J. D. McCarthy and A. Papadopoulos, {\it Tilings of surfaces: the complex of domains}, monograph in preparation

\bibitem{mv2} J. D. McCarthy and W. R. Vautaw, Automorphisms of the complex of separating curves in genus $3$, to appear, preliminary report - Braids, links, and mapping class groups, a conference in honor of Professor Joan Birman, Columbia University, New York, March 17, 2005.

\bibitem{munkres} J. R. Munkres, {\it Elements of algebraic topology}, Addison-Wesley Publishing Company, Menlo Park, Calif., 1984, 454 pp., ISBN 0-201-04586-9, MR85m:55001

\bibitem{shackleton} K. J. Shackleton, Combinatorial rigidity in curve complexes and mapping class groups, LANL ArXiv, math.GT/0503199

\end{thebibliography}
\end{document}